\documentclass[11pt]{article}

\usepackage{amsmath,amsthm,epsfig}

\def\ftoday{le \space\number\day \space\ifcase\month\or

  janvier\or f\'evrier\or mars\or avril\or mai\or juin\or

  juillet\or ao\^ut\or septembre\or octobre\or novembre\or d\'ecembre\fi

  \space\number\year}

\def\real{I\kern-0.20em R}

\def\integer{I\kern-0.20em N}

\def\relative{{\rm \rlap Z\kern 2.2pt Z}}

\def\cc{\kern-.25em{\c c}}

\def\bc{\begin{center}}

\def\ec{\end{center}}

\def\=def{\stackrel{{\rm def}}{=}}

\newcommand\vf[1]{{\cal X}_{#1}}

\newcommand\vfo[2]{{\cal X}_{#1}^{#2}}

\newcommand\fvfo[2]{\widehat{\cal X}_{#1}^{#2}}

\newcommand\lie[1]{{\frak #1}}

\newcounter{indconst}

\newcounter{auxconst}

\def\bit{\begin{itemize}}

\def\eit{\end{itemize}}

\def\ben{\begin{enumerate}}

\def\een{\end{enumerate}}

\def\bde{\begin{description}}

\def\ede{\end{description}}

\def\beq{\begin{equation}}

\def\eeq{\end{equation}}

\def\bfi{\begin{figure}[hbt] \begin{center}}

\def\efi{\end{center} \end{figure}}

\def\bce{\begin{center}}

\def\ece{\end{center}}

\newtheorem {theo} {Th\'eor\`eme}[section]

\newtheorem {coro} {Corollaire} [section]

\newtheorem {lemm} {Lemme}[section]

\newtheorem {prop} {Proposition}[section]

\newtheorem {defi} {D\'efinition}[section]

\newtheorem {rem}{Remarque} [section]

\newtheorem {theos} {Theorem}[subsection]

\newtheorem {ex}[theos] {Exemple}

\input amssym.def

\input amssym

\numberwithin{equation}{section}

\def\abstractname{R\'{e}sum\'{e}}

\begin{document}
\title{Sur les structures de Poisson singuli\`eres}
\author{Laurent Stolovitch \thanks{CNRS UMR 5580, Laboratoire Emile Picard,
Universite Paul Sabatier, 118 route de Narbonne,
31062 Toulouse cedex 4, France. Courriel : {\tt stolo@picard.ups-tlse.fr}}}
\date{\ftoday}
\maketitle
\begin{abstract}
Nous nous int\'eressons aux structures de Poisson analytiques 
singuli\`eres en un point et de partie lin\'eaire non-nulle en ce point. En utilisant des travaux r\'ecents de l'auteur sur la
normalisation holomorphe de familles commutatives de champs de vecteurs holomorphes singuliers, nous donnons des r\'esultats de normalisation holomorphe de certaines structures de Poisson.
\end{abstract}
\def\abstractname{Abstract}
\begin{abstract}
We are interested in analytic singular Poisson structures with a non zero linear part
at the singularity. Using recent work of the author about holomorphic
normalization of commutative familly of singular vector fields, we obtain
results about normalization of holomorphic Poisson structures.
\end{abstract}
\section{Introduction}

Soit M une vari\'et\'e analytique de dimension N. Une {\bf structure de Poisson} est la donn\'ee d'un crochet $\{.,.\}$ qui assigne, \`a un couple $(f,g)$ de germes de fonctions holomorphes en un point $x$ de $M$, un germe $\{f,g\}$ de fonction holomorphe en $x$ v\'erifiant les propri\'et\'es suivantes~:
\begin{itemize}
\item $\{.,.\}$ est bilin\'eaire et antisym\'etrique,
\item $\{f,\{ g, h\}\}+\{g,\{ h, f\}\}+\{h,\{ f, g\}\}=0$ (identit\'e de Jacobi),
\item $\{f,gh\}=\{f,g\}h+\{f,h\}g$ (identit\'e de Leibniz).
\end{itemize}

Il revient au m\^eme de d\'efinir un champ de $2$-vecteurs que l'on peut \'ecrire, dans un syst\`eme de coordonn\'ees locales,
$$
P=\frac{1}{2}\sum_{1\leq i,j\leq N}P_{i,j}(x)\frac{\partial}{\partial x_i}\wedge \frac{\partial}{\partial x_j}=\sum_{1\leq i<j\leq N}P_{i,j}(x)\frac{\partial}{\partial x_i}\wedge \frac{\partial}{\partial x_j}\quad\text{avec }\;P_{i,j}=-P_{j,i}
$$
et qui v\'erifie l'identit\'e de Jacobi
$$
\sum_{1\leq l\leq N}\left(P_{i,l}\frac{\partial P_{j,k}}{\partial x_l}+P_{j,l}\frac{\partial P_{k,i}}{\partial x_l}+P_{k,l}\frac{\partial P_{i,j}}{\partial x_l}\right)=0
$$
pour $1\leq i,j,k\leq N$. On d\'efinit alors le crochet de Poisson par 
$$
\{f,g\}:=<P,df\wedge dg>= \sum_{1\leq i<j\leq N}P_{i,j}(x)\left(\frac{\partial f}{\partial x_i} \frac{\partial g}{\partial x_j}-\frac{\partial g}{\partial x_i} \frac{\partial f}{\partial x_j}\right).
$$
On a alors $\{x_i,x_j\}=P_{i,j}$. Le hamiltonien associ\'e \`a un germe de fonction $f$ relativement \`a $P$ est le germe de champ de vecteurs d\'efini par
$$
X_f=\sum_{i=1}^N\{x_i,f\}\frac{\partial }{\partial x_i}.
$$
Soit $g$ un autre germe de fonction holomorphe. La d\'eriv\'ee de Lie de $g$ le long de $X_f$, $X_f(g)$, est \'egale \`a $\{f,g\}$.

D'apr\`es un r\'esultat d'A. Weinstein \cite{weinstein1, weinstein-cannas}, on peut trouver un bon syst\`eme de coordonn\'ees holomorphes $(p_1,\ldots, p_m,q_1,\ldots, q_m, x_1,\ldots, x_r)$ au voisinage d'un point $p$ dans lequel on ait $p_i(p)=0, q_i(p)=0, x_j(p)=0$ et 
$$
P=\sum_{1\leq i\leq m}\frac{\partial}{\partial p_i}\wedge \frac{\partial}{\partial q_i}+\sum_{1\leq i<j\leq r}P_{i,j}(x)\frac{\partial}{\partial x_i}\wedge \frac{\partial}{\partial x_j}
$$
avec $P_{i,j}(0)=0$ et $2m+r=N$.

L'\'etude de la classification holomorphe locale d'un crochet de Poisson
(c'est-\`a-dire les classes d'\'equivalence par conjugaison par des
diff\'eomorphismes locaux holomorphes) revient donc \`a celle d'un crochet de
Poisson nul en un point. Pour le cas de la  dimension deux, on pourra
consulter \cite{arn2-2nd}[Appendix 14]

{\bf Dans la suite, nous supposerons que $M$ est un voisinage de l'origine dans $\Bbb C^N$ et que $P$ s'annule en ce point.}

Posons
$$
c_{i,j}^k=\frac{\partial P_{i,j}}{\partial x_k}(0)\quad\quad 1\leq i,j,k\leq N.
$$
La partie lin\'eaire (ou $1$-jet en $0$) de $P$ d\'efinit une structure d'alg\`ebre de Lie sur l'espace cotangent en $0$, ${\lie g}$, de la mani\`ere suivante :
$$
[x_i,x_j]=\sum_{i=1}^Nc_{i,j}^kx_k.
$$
L'identit\'e de Jacobi se montre en prenant le $1$-jet, en $0$, de l'identit\'e de Jacobi satisfaite par $P$.

Lorsque que cette alg\`ebre est semi-simple, J. Conn \cite{conn1, conn2} a d\'emontr\'e qu'un tel crochet est holomorphiquement lin\'earisable. Cela signifie qu'il existe un syst\`eme de coordonn\'ees holomorphes, $y=\phi(x)$, dans lequel on a
$$
\phi_*P(y)=\sum_{1\leq i<j\leq n}\left(\sum_{k=1}^Nc^k_{i,j}y_k\right)\frac{\partial}{\partial y_i}\wedge \frac{\partial}{\partial y_j}.
$$

R\'ecemment, Nguyen Tien Zung \cite{zung1} a donn\'e une d\'ecomposition de
Levi holomorphe d'un crochet de Poisson singulier, montrant ainsi que "la difficult\'e r\'eside" dans le radical de l'alg\`ebre.

La majeure partie des travaux portant sur l'\'etude locale d'une structure de
Poisson concerne les cas de lin\'earisation d'une telle structure (par
exemple \cite{dufourmolinier, dufourzungcras}). En particulier, B. Abbaci \cite{abbaci} a am\'elior\'e, dans une certaine situation, certains r\'esultats de J.-P. Dufour.

Dans cet article, nous nous proposons non pas de lin\'eariser mais de normaliser holomorphiquement certaine structures de Poisson dont l'alg\`ebre de Lie associ\'ee est un produit semi-direct $\Bbb C^p\ltimes \Bbb C^n$. 

Notre travail s'inspire des articles de J.-P. Dufour \cite{dufour2} et de
Dufour-Zhitomirskii \cite{dufour3}.

\subsection{Notations}

Soit $n$ un entier non-nul.
On notera~:
\begin{itemize}
\item $\vfo n k$ (resp. $\fvfo n k$) l'espace des germes en $0\in \Bbb C^n$ de
  champs de vecteurs holomorphes (resp. formels) d'ordre $\geq k$ en $0$;
\item $\vf {n,p}:=\wedge^p \vf n$ l'espace des germes en $0\in \Bbb C^n$ de champs de tenseurs contravariants holomorphes anti-sym\'etriques de type $(p,0)$;
\item ${\cal O}_n$ (resp. $\widehat{\cal O}_n$) l'anneau des germes en $0\in
  \Bbb C^n$ de fonctions holomorphes (resp. formelles);
\item ${\cal M}_n$ l'id\'eal maximal de ${\cal O}_n$;
\item ${\cal O}_n(U)$ (resp. ${\cal O}_n(K)$) l'anneau des fonctions
  holomorphes sur l'ouvert $U$ (resp. au voisinage du compact $K$) de $\Bbb C^n$;
\item si $X\in \fvfo n 1$ et $k\in \Bbb N^*$, $J^k(X)$ d\'esigne le polyn\^ome
  de Taylor \`a l'ordre  $k$ en $0$ de $X$.
\end{itemize}

Soient $Q=(q_1,\ldots,q_n)\in \Bbb N^n$ et $\lambda=(\lambda_1,\ldots,\lambda_n)\in \Bbb N^n$. On posera $|Q|=q_1+\ldots+q_n$ et 
$(Q,\lambda)=q_1\lambda_1+\ldots+ q_n\lambda_n$. Soit $k\in \Bbb N$, on notera 
$\Bbb N^n_k$, l'ensemble des multiindices $Q\in \Bbb N^n$ tels que $|Q|\geq k$.

Soit $K$ un compact de $\Bbb C^p$. Un \'el\'ement $f$ de ${\cal O}_p(K)\otimes
{\cal O}_n$ d\'efini une fonction holomorphe que l'on d\'eveloppera "le long"
 de $K\times\{0\}$ de la
mani\`ere suivante :
$$
f=\sum_{Q\in N^n}f_Q(x'')(x')^Q
$$
o\`u $x'$ (resp. $x''$) d\'esigne les coordonn\'ees dans $\Bbb C^n$ (resp.
$\Bbb C^p$) et $f_Q\in {\cal O}_p(V)$, $V$ \'etant un voisinage de $K$
ind\'ependant de $Q$.

\subsection{Formes normales de champs de vecteurs}\label{forme-normale}

Un germe de champ de vecteurs $X$ de $(\Bbb C^n,0)$, nul en ce point mais de
partie lin\'eaire $L$ non nulle,  sera dit normalis\'e, ou sous {\bf forme
normale}, s'il commute avec sa partie lin\'eaire, i.e. $[L, X]=0$. 
On rappelle le th\'eor\`eme de Poincar\'e-Dulac \cite{Arn2}
\begin{theo}
Soient $X=S+R$ un champ de vecteurs formel de $(\Bbb C^n,0)$ s'annulant \`a l'origine; $S$ sa partie
lin\'eaire suppos\'ee semi-simple et $R$ 
un champ de vecteurs non-lin\'eaire. Il existe un diff\'eomorphisme formel
$\hat\Phi$ de $(\Bbb C^n,0)$, tangent 
\`a l'identit\'e en $0$ tel que $\hat\Phi_*X=S+N$ o\`u $N$ est un champ de
vecteurs formel tel que $[S,N]=0$. 
\end{theo}
On dit alors que $\hat\Phi_*X(y):=D\Phi(\Phi^{-1}(y))X(\Phi^{-1}(y))$ est une forme normale formelle de $X$.

Soit $X=\{X_1,\ldots, X_p\}$ une famille commutative de champs de vecteurs
holomorphes au voisinage de l'origine dans $\Bbb C^n$. On les suppose nuls en
ce point et on suppose que leurs parties lin\'eaires sont diagonales et
lin\'eairement ind\'ependantes sur $\Bbb C$. On posera $J^1(X_i)=S_i$. On peut d\'efinir une
notion de forme normale formelle de la famille $\{X_i\}$ relativement \`a la
famille de parties lin\'eaires $S=\{S_i\}$ \cite{stolo-ihes}[propostion
4.2.1]. Nous dirons que le famille $X$ est formellement conjugu\'ee \`a la
famille $\{Y_i\}_{1\leq i\leq p}$ s'il existe un diff\'eomorphisme formel
$\hat\Phi$ de $(\Bbb C^n,0)$ fixant $0$ et tangent \`a l'identit\'e en ce
point et tel que $\hat\Phi_*X_i=Y_i$ pour $1\leq i\leq p$. 

Soit ${\lie g}$ une alg\`ebre de Lie commutative de dimension $p$ sur $\Bbb
C$. Les familles $S$ et $X$ d\'efinissent des morphismes de Lie de ${\lie g}$
dans l'alg\`ebre de Lie des germes de champs de vecteurs holomorphes nuls \`a
l'origine de $\Bbb C^n$ de la mani\`ere suivante :  soit $\{g_1,\ldots, g_p\}$
une base de ${\lie g}$. On pose alors $X(g_i)=X_i$ et $S(g_i)=S_i$.
Soient $\lambda_1,\ldots,\lambda_n$ des formes lin\'eaires complexes sur $\lie
g$ tel que le morphisme de Lie $S$ de $\lie g$ 
dans l'alg\`ebre de Lie des champs de vecteurs lin\'eaires de $\Bbb C^n$
d\'efini par 
$S(g)=\sum_{i=1}^n\lambda_i(g)x_i\frac{\partial}{\partial x_i}$ soit
injectif. Pour tout $Q=(q_1,\ldots,q_n)\in \Bbb N^n$ et $1\leq i\leq n$, 
on d\'efinit le {\bf poids} $\alpha_{Q,i}(S)$ de $S$ comme \'etant la forme
lin\'eaire $\sum_{j=1}^n{q_j\lambda_j(g)}-\lambda_i(g)$. L'espace de poids
associ\'e \`a un poids $\alpha$ est l'espace
$$
\vfo {n,\alpha} {1}=\left\{p\in \vfo n 1\,|\,\forall g \in \lie g,\,\,[S(g),p]=\alpha(g)p\right\}.
$$

D\'efinissons les espaces de poids nul
\begin{eqnarray*}
\widehat{\cal O}_n^S := \left\{f\in\widehat {\cal O}_n\;|\; S_i(f)=0\; 1\leq i\leq p \right\},\\
(\fvfo n 1)^S:= \left\{X\in \fvfo n 1\;|\; [S_i,X]=0\; 1\leq i\leq p \right\},\\
\end{eqnarray*}
o\`u $S_i(f)$ d\'esigne la d\'eriv\'ee de Lie de $f$ le long de $S_i$.

Si l'anneau ${\cal O}_n^S$ n'est pas r\'eduit au corps des nombres complexes, on sait \cite{stolo-ihes}[proposition 5.3.2] qu'il existe un nombre fini de mon\^omes $u_1:=x^{R_1},\ldots, u_t:=x^{R_t}$, $R_i\in \Bbb N^n$, tels que
$\widehat{\cal O}_n^S=\Bbb C[[u_1,\ldots, u_t]]$. D'autre part, $(\fvfo n 1)^S$ est un module de type fini sur $\widehat{\cal O}_n^S$.

\subsection{Forme normale de crochet de Poisson}

Dans l'article de J. Conn \cite{conn1}, la notion de forme normale de
structure de Poisson est esquiss\'ee m\^eme si elle n'est utilis\'ee que pour
la lin\'eariser. Elle a \'et\'e ensuite d\'evelopp\'ee par O.V. Lychagina
\cite{lychagina} (l'auteur remercie J.-P. Dufour d'avoir port\'e cet article \`a sa
connaissance).
La structure de Poisson lin\'earis\'ee d\'efinit une repr\'esentation de
l'alg\`ebre de Lie associ\'ee ${\lie g}$ dans chaque espace de polyn\^omes
homog\`enes sur ${\lie g}$. Elle est d\'efinie par la d\'eriv\'ee de Lie le
long du hamiltonien associ\'e \`a un \'element de ${\lie g}$, c'est -\`a-dire
$$
\rho(x_i)(f)=X_{x_i}(f)\quad\quad 1\leq i\leq N.
$$
On d\'efinit alors le complexe de Chevalley-Kozsul et ses groupes de
cohomologie $H^i({\lie g}, S^k(\lie g ))=Z^i({\lie g}, S^k(\lie g))/B^i({\lie
  g}, S^k(\lie g))$ \cite{serre-lie}, $S^k(\lie g)$ \'etant la $k$-i\`eme puissance sym\'etrique de ${\lie g}$ \cite{stolo-ihes}[section 4.1]. 
On a alors le r\'esultat suivant
\begin{prop}\cite{conn1, lychagina}
Soit $P$ une structure de Poisspin.texson nulle en $0\in \Bbb C^N$ et de $1$-jet non-nul en ce point. Pour tout entier $k\geq 2$, soit $V_k$ un suppl\'ementaire de $B^2({\lie g}, S^k(\lie g))$ dans $Z^2({\lie g}, S^k(\lie g))$. Il existe alors un diff\'eomorphisme formel $\hat \Phi$ de $(\Bbb C^N,0)$ nul en $0$ et tangent \`a l'identit\'e en ce point tel que 
$$
\hat\Phi_*P-{\cal L}\in\oplus_{k\geq 2}V_k.
$$
o\`u ${\cal L}$ d\'esigne la partie lin\'eaire de $P$. 
\end{prop}
Nous dirons alors que $\hat\Phi_*P$ est {\bf normalis\'e}.

On d\'efinit le {\bf crochet de Schouten-Nijenhuis} $[.,.]: \vf {N,p}\times \vf {N,q}\rightarrow \vf {N,p+q-1}$, $\Bbb C$-bilin\'eaire, qui est une extension de la d\'eriv\'ee de Lie et qui v\'erifie
\begin{enumerate}
\item $[P,Q]=(-1)^{pq}[Q,P]$,\\
\item $[P, Q\wedge R]=[P,Q]\wedge R + (-1)^{pq+q}Q\wedge [P,R]$,
\item $(-1)^{p(r-1)}[P,[Q,R]]+ (-1)^{q(p-1)}[Q,[R,P]] +(-1)^{r(q-1)}[R, [P,Q]]=0$.
\end{enumerate}
Nous renvoyons le lecteur aux ouvrages \cite{vaisman, weinstein-cannas, schouten} pour de plus amples renseignements et
la bibliographie. On a les propri\'et\'es suivantes. Un champ de $2$-tenseur $P$ d\'efinit une structure de Poisson
si et seulement si il v\'erifie $[P,P]=0$. De plus, si $X_f$ est un champ
hamiltonien associ\'e \`a une structure de Poisson $P$ alors $[P, X_f]=0$.

Le morphisme $S$ d\'efinit pr\'ec\'edement induit, par le crochet de Schouten-Nijenhuis, une
repr\'esentation de ${\lie g}$ dans l'espace des germes de champs de
bi-vecteurs holomorphes et nuls \`a l'origine gr\^ace \`a l'identit\'e de
"Jacobi" ci-dessus. 
\begin{lemm}
\begin{enumerate}
\item Les poids de $S$ dans $\vfo {n,2} 2$ sont de la forme
  $(Q,\lambda(g))-\lambda_i(g)-\lambda_j(g)$, $Q\in \Bbb N^n_2$, $1\leq i,j\leq
  n$.
\item Soient $\alpha $ et $\beta$ deux poids de $S$ dans $\vfo n 1$.
Soient $X\in \vfo {n,\alpha} 1$ et $Y\in \vfo {n,\beta} 1$. Alors, $X\wedge Y$
est de poids $\alpha+\beta$. 
\end{enumerate}
\end{lemm}\label{poids-poisson}
On a 
\begin{eqnarray*}
\left [S(g), x^Q\frac{\partial}{\partial x_i}\wedge\frac{\partial}{\partial
    x_j}\right] & = & \left[S(g), x^Q\frac{\partial}{\partial x_i}\right]\wedge\frac{\partial}{\partial
    x_j}+ x^Q\frac{\partial}{\partial x_i}\wedge \left[S(g),\frac{\partial}{\partial
    x_j}\right]\\
& = & \left((Q,\lambda(g))-\lambda_i(g)-\lambda_j(g)\right) x^Q\frac{\partial}{\partial x_i}\wedge\frac{\partial}{\partial
    x_j}.
\end{eqnarray*}
En effet, on a 
\begin{eqnarray*}
[S(g),X\wedge Y] & = & [S(g),X]\wedge Y + X\wedge [S(g),Y]\\
& = & (\alpha(g)+\beta(g))X\wedge Y.
\end{eqnarray*}
Les espaces de poids sont en somme directes (en fait, ce sont les espaces de
poids dans les espaces de champs de bi-vecteurs homog\`enes qui sont en somme
directe \cite{Lie78-bourbaki}).


\section{Cadre de travail et premi\`eres r\'eductions}\label{red}

Nous supposerons que l'alg\`ebre de Lie ${\lie g}=\Bbb C^p\ltimes\Bbb C^n$ est le produit semi-direct de $\Bbb C^p$ par $\Bbb C^n$ (la suite de morphismes d'alg\`ebres de Lie $0\rightarrow \Bbb C^n\rightarrow \lie g\rightarrow \Bbb C^p\rightarrow 0$ est exacte) d\'efini par 
$$
\left\{
\begin{array}{lcc}
\;[x_i,x_j]= 0 & &1\leq i,j\leq n\\ 
\;[x_i,x_j]=0 & & n+1\leq i,j\leq N\\ 
\;[x_i,x_j]=\lambda_{j,i}x_i & & 1\leq i\leq n<j\leq N
\end{array}\right.
$$
o\`u l'on a pos\'e $N=n+p$.
On supposera que les champs de vecteurs 
$$
S_j:=\sum_{i=1}^n\lambda_{j,i}x_i\frac{\partial}{\partial x_i}\quad 1\leq j\leq p
$$ 
sont lin\'eairement ind\'ependants sur $\Bbb C$. On notera $S$ la famille des $S_j$ et $\lambda^i=(\lambda_{i,1},\ldots, \lambda_{i,n})$.

{\bf Dans la suite nous ferons les hypoth\`eses suivantes (H)}~:
\begin{enumerate}
\item un des $S_j$ a des valeurs propres distinctes, i.e. $x_j\partial/\partial
  x_i\not\in (\fvfo n 1)^S$ si $i\neq j$;
\item un des $S_j$ n'a pas de valeur propre nulle, i.e. pour $1\leq i\leq n$, $x_i\not\in \widehat{\cal O}_n^S$;
\item un des $S_j$ n'a pas deux valeurs propres de somme nulle, i.e. pour $1\leq i,j\leq n$, $x_ix_j\not\in \widehat{\cal O}_n^S$;
\item un des $S_j$ est tel qu'aucune de ses valeurs propres n'est la somme de
  deux autres, i.e. pour $1\leq i,j,k\leq n$, $x_ix_j\partial/\partial
  x_k\not\in (\fvfo n 1)^S$;
\item $\{x_i,x_j\}=0$ pour $n+1\leq i,j\leq N$.
\end{enumerate} 

Pour tout $1\leq j\leq p$, on notera $X_j$ le hamiltonien asscoci\'e \`a $x_{n+j}$. Par hypoth\`ese, on a 
$$
X_j=\sum_{i=1}^n\{x_i,x_{n+j}\}\frac{\partial}{\partial x_i}.
$$
On a alors
$$
P=\sum_{k=1}^pX_k(x)\wedge \frac{\partial}{\partial x_{n+k}} + \sum_{1\leq i<j\leq n} \{x_i,x_{j}\}\frac{\partial}{\partial x_i}\wedge \frac{\partial}{\partial x_{j}}.
$$
Soit $f$ un germe de fonction et $1\leq j,k\leq p$. L'identit\'e de Jacobi
$$
\{x_{n+j},\{x_{n+k},f\}\}=\{x_{n+k},\{x_{n+j},f\}\},
$$
montre que les champs $X_j$ commutent deux-\`a-deux. Les champs $X_j$ peuvent
\^etre consid\'er\'es comme des germes de champs de $\Bbb C^n$ d\'ependant
holomorphiquement de $p$ param\`etres $x'':=(x_{n+1},\ldots, x_N)$. On notera
${\cal O}_p=\Bbb C\{x_{n+1},\ldots, x_N\}$ (resp. ${\cal O}_n=\Bbb
C\{x_{1},\ldots, x_n\}$) l'anneau des germes de fonctions holomorphes de
$(\Bbb C^p,0)$ (resp. $(\Bbb C^n,0)$). On posera $x':=(x_{1},\ldots, x_n)$ et
on notera ${\cal M}$ l'id\'eal maximal de ${\cal O}_n$. Soit $f\in {\cal
  O}_p\otimes {\cal O}_n$, $f^k(x)$ est la composante homog\`ene de degr\'e
$k$ en $x'=0$ de $f$.
 On peut alors \'ecrire pour $1\leq j\leq p$,
$$
X_j(x)=\sum_{i=1}^nb_{j,i}(x'')\frac{\partial}{\partial x_i}+\sum_{i,k=1}^n b_{j,i,k}(x'')x_k\frac{\partial}{\partial x_i} \mod{ {\cal O}_p \otimes {\cal M}^2}.
$$
Par hypoth\`ese, la diff\'erentielle de $X_j^0$ en $0$ est nulle et 
$$
X_j^1(x',0)=S_j=\sum_{i,k=1}^n b_{j,i,k}(0)x_k\frac{\partial}{\partial x_i}.
$$

\begin{lemm}

Quitte \`a faire un changement de variables holomorphe de $(\Bbb C^N,0)$ tangent \`a l'identit\'e en $0$, on peut supposer que $X_j(0,x'')=0$ et $DX_j(0,x'')=\text{diag}(\lambda_{j,1}(x''),\ldots, \lambda_{j,n}(x''))$ pour $1\leq j\leq p$.

\end{lemm}

\begin{proof}

Par hypoth\`ese, on peut supposer que la matrice $(b_{1,i,k}(0))_{1\leq i,k\leq n}$ est inversible; en appliquant le th\'eor\`eme des fonctions implicites, on obtient un germe d'application holomorphe $g:(\Bbb C^p,0)\rightarrow (\Bbb C^n,0)$ telle que $g(0)=0$ et $X_1(g(x''),x'')=0$. Le champ $X_1$ est donc tangent au graphe ${\cal G}$ de $g$; puisque les $X_i$ commutent \`a $X_1$, ils sont aussi tangents \`a ${\cal G}$. Les $X_j$ n'ont aucune composante non nulle le long des $\frac{\partial }{\partial x_{n+k}}$, $1\leq k\leq p$; ils doivent donc \^etre nuls sur ${\cal G}$. Faisons le changement de variables holomorphe $(X',X'')=(x'-g(x''),x'')$ fixant l'origine et tangent \`a l'identit\'e en ce point. On a alors $X_j(0,X'')=0$ et $DX_j(0)=(b_{j,i,k}(0))_{j\leq i,k\leq n}$. 

Puisque les valeurs propres de l'un des $DX_i(0)$ sont distinctes, on peut
diagonaliser $DX_i(0,x'')$ holomorphiquement sur un voisinage $U$ de $0$ dans
$\Bbb C^p$. Les matrices $DX_j(0,x'')=(b_{j,i,k}(x''))_{j\leq i,k\leq n}$ sont
alors diagonalis\'ees et leurs coefficients sont holomorphes dans $U$. Ce qui donne le r\'esultat.

\end{proof}
%
%

\begin{lemm}\label{ordre2}

Sous les hypoth\`eses pr\'ec\'edentes, on a $\{x_i,x_{j}\}\in {\cal O}_p \otimes {\cal M}^2$ pour $1\leq i,j\leq n$.

\end{lemm}

\begin{proof}

Posons 
$$
\tilde P = \sum_{1\leq i<j\leq n} \{x_i,x_{j}\}\frac{\partial}{\partial x_i}\wedge \frac{\partial}{\partial x_{j}}
$$
Pour all\'eger l'\'ecriture, on posera $g_{i,j}= \{x_i,x_{j}\}$. Puisque $X_i$, $1\leq i\leq p$, est un hamiltonien par rapport \`a $P$, le crochet de Schouten-Nijenhuis $[P,X_i]$ est nul. On a donc 

\begin{eqnarray}\label{pxi}
 [\tilde P, X_i] & = &\sum_{k=1}^p \frac{\partial X_i}{\partial x_{n+k}}\wedge
X_k(x)
\end{eqnarray}

Or, on a 

\begin{eqnarray*}
[\tilde P, X_i] & = & \sum_{1\leq r< s\leq
   n}\left[g_{r,s}(x)\frac{\partial}{\partial x_r}\wedge
   \frac{\partial}{\partial x_{s}}, X_i\right]\\
& = & \sum_{1\leq r< s\leq n}X_i(g_{r,s})\frac{\partial}{\partial x_r}\wedge \frac{\partial}{\partial x_{s}}-g_{r,s}\left[ X_i, \frac{\partial}{\partial x_r}\wedge \frac{\partial}{\partial x_{s}}\right].
\end{eqnarray*}

Or, un calcul simple montre que l'on a

\begin{eqnarray*}
\left[ X_i, \frac{\partial}{\partial x_r}\wedge \frac{\partial}{\partial x_{s}}\right] &= & \sum_{k=1}^n \left(\frac{\partial X_{i,k}}{\partial x_r}\frac{\partial}{\partial x_k}\wedge \frac{\partial}{\partial x_{s}}+ \frac{\partial X_{i,k}}{\partial x_s}\frac{\partial}{\partial x_r}\wedge \frac{\partial}{\partial x_{k}}\right).
\end{eqnarray*}

On obtient alors 

\begin{eqnarray}\label{tildepxi}
[\tilde P, X_i] & = & \sum_{1\leq r< s\leq n}\left(X_i(g_{r,s})-g_{r,s}\left(\frac{\partial X_{i,r}}{\partial x_r}+\frac{\partial X_{i,s}}{\partial x_s}\right)\right)\frac{\partial}{\partial x_r}\wedge \frac{\partial}{\partial x_{s}}.
\end{eqnarray}

D'apr\`es le lemme pr\'ec\'edent, la restriction \`a $x'=0$ du membre de droite de l'\'equation $(\ref{pxi})$ est nulle. Or, la restriction \`a $x'=0$ du membre de droite de l'\'equation $(\ref{tildepxi})$ est \'egale \`a

$$
-\sum_{1\leq r< s\leq n}g_{r,s}^0(x'')(\lambda_{i,r}(x'')+\lambda_{i,s}(x''))\frac{\partial}{\partial x_r}\wedge \frac{\partial}{\partial x_{s}}.
$$

Donc, pour $1\leq i\leq p$ et pour $1\leq r<s\leq n$, on a 

$$
g_{r,s}^0(x'')(\lambda_{i,r}(x'')+\lambda_{i,s}(x''))=0
$$

Par hypoth\`ese, pour $1\leq r<s\leq n$, il existe un entier $1\leq i\leq
p$ tel que $\lambda_{i,r}(0)+\lambda_{i,s}(0)\neq 0$. Donc,
$\lambda_{i,r}(x'')+\lambda_{i,s}(x'')$ est non-nul sur un voisinage de l'origine. Par cons\'equent, la fonction holomorphe
$g_{r,s}^0(x'')$ est nulle sur le connexe $U$.

Le 1-jet en $x'=0$ du membre de droite de l'\'equation $(\ref{pxi})$ est nul. Celui du membre de droite de l'equation $(\ref{tildepxi})$ est \'egal \`a 
$$
\sum_{1\leq r< s\leq n}\left(J^1(X_i)(g_{r,s}^1)-g^1_{r,s}(\lambda_{i,r}(x'')+\lambda_{i,s}(x''))\right)\frac{\partial}{\partial x_r}\wedge \frac{\partial}{\partial x_{s}}.
$$
En posant $g^1_{r,s}=\sum_{k=1}^ng^1_{r,s,k}(x'')x_k$, on obtient pour $1\leq r,s,k\leq n$,
$$
g^1_{r,s,k}(x'')(\lambda_{i,k}(x'')-\lambda_{i,r}(x'')-\lambda_{i,s}(x''))=0.
$$
Par hypoth\`ese, pour chaque triplet $(r,s,k)$, il existe un entier $1\leq
i\leq n$ tel que $\lambda_{i,k}(0)-\lambda_{i,r}(0)-\lambda_{i,s}(0)\neq
0$. Les fonctions $g^1_{r,s,k}(x'')$ sont donc toutes nulles.

Par cons\'equent, nous avons montr\'e que pour $1\leq i,j\leq n$,
$$
\{x_i,x_j\}\in {\cal O}_p \otimes {\cal M}^2.
$$
\end{proof}

\begin{lemm}

Il existe un changement de variables holomorphe $(X',X'')=\Psi(x)=(x',\psi(x''))$ de $(\Bbb C^N,0)$ tel que, pour $1\leq i\leq p$, la partie lin\'eaire de $\Psi_*X_i$ en $0$ est \'egale \`a $S_i$ (les valeurs propres sont ind\'ependantes de $x''$).

\end{lemm}

\begin{proof}

Prenons le $2$-jet de l'\'egalit\'e $(\ref{pxi})$. Le $2$-jet de $(\ref{tildepxi})$ est \'egal \`a 

$$
\sum_{1\leq r< s\leq n}\left(J^1(X_i)(g_{r,s}^2)-g^2_{r,s}(\lambda_{i,r}(x'')+\lambda_{i,s}(x''))\right)\frac{\partial}{\partial x_r}\wedge \frac{\partial}{\partial x_{s}},
$$

c'est-\`a-dire 

$$
\sum_{1\leq r< s\leq n}\left(\sum_{Q\in \Bbb N^n, |Q|=2}\left((Q,\lambda_i(x''))-\lambda_{i,r}(x'')-\lambda_{i,s}(x'')\right)g_{r,s,Q}^2(x'')(x')^Q\right)\frac{\partial}{\partial x_r}\wedge \frac{\partial}{\partial x_{s}}
$$

o\`u l'on a \'ecrit 

$$
g_{r,s}^2(x)=\sum_{Q\in \Bbb N^n, |Q|=2}g_{r,s,Q}^2(x'')(x')^Q.
$$

D'autre part, le $2$-jet du membre de droite de $(\ref{pxi})$ est \'egal \`a

\begin{eqnarray*}
\sum_{k=1}^p \frac{\partial J^1(X_i)}{\partial x_{n+k}}\wedge J^1(X_k)(x) & = & \sum_{1\leq r< s\leq n}\left(\sum_{k=1}^p\frac{\partial \lambda_{i,r}}{\partial x_{n+k}}\lambda_{k,s}-\frac{\partial \lambda_{i,s}}{\partial x_{n+k}}\lambda_{k,r}\right)x_rx_s\frac{\partial}{\partial x_r}\wedge \frac{\partial}{\partial x_{s}}\\
\end{eqnarray*}

En comparant les coefficients de $x_rx_s\frac{\partial}{\partial x_r}\wedge \frac{\partial}{\partial x_{s}}$ dans les deux expressions, on obtient 

$$
\sum_{k=1}^p\left(\frac{\partial \lambda_{i,r}}{\partial x_{n+k}}\lambda_{k,s}-\frac{\partial \lambda_{i,s}}{\partial x_{n+k}}\lambda_{k,r}\right)=0\quad 1\leq r< s\leq n.
$$

En posant $\Lambda_m:=\sum_{i=1}^p\lambda_{i,m}(x'')\frac{\partial }{\partial x_{n+i}}$ pour $1\leq m\leq n$, les \'equations pr\'ec\'edentes s'\'ecrivent $[\Lambda_r,\Lambda_s]=0$ pour $1\leq r,s\leq n$. Or, par hypoth\`ese, la famille $\{\Lambda_r(0)\}_{1\leq r\leq n}$ est de rang $p$. On peut supposer que la famille $\{\Lambda_r(0)\}_{1\leq r\leq p}$ est ind\'ependante sur $\Bbb C$.

On peut donc redresser simultan\'ement et holomorphiquement la famille $\{\Lambda_r(x'')\}_{1\leq r\leq p}$~: il existe un changement de variables holomorphe $X''=\psi(x'')$ tel que $\psi_*\Lambda_r(X'')=\Lambda_r(0)$ pour $1\leq r\leq p$. Soit $p+1\leq l\leq n$; pour $1\leq r\leq p$, on a 

\begin{eqnarray*}
0=[\psi_*\Lambda_r, \psi_*\Lambda_l] & =: &  \left[ \Lambda_r(0), \sum_{j=1}^p\mu_{j,l}(X'')\frac{\partial}{\partial X_{j+n}}\right]\\
& = & \sum_{j=1}^p\Lambda_r(0)(\mu_{j,l}(X''))\frac{\partial}{\partial X_{j+n}}.
\end{eqnarray*}
o\`u les $\mu_{j,l}$ sont des germes de fonctions holomorphes dans $(\Bbb C^p,0)$.
De par l'ind\'ependance des $p$ premiers $\Lambda_r(0)$, l'\'equation pr\'ec\'edente se ram\`ene \`a ${\partial \mu_{j,l}}/{\partial X_{k+n}}=0$ pour $1\leq k\leq p$. Donc, $\psi_*\Lambda_l(X'')=\Lambda_l(0)$.

\end{proof}


\section{R\'esultats principaux}

Soit $\{.,.\}$ un crochet de Poisson holomorphe de $(\Bbb C^N,0)$, nul en $0$ et v\'erifiant les hypoth\`eses $(H)$ d\'efinie dans la section $\ref{red}$. D'apr\`es les r\'esultats pr\'ec\'edents, quitte \`a
faire un changement de coordonn\'ees holomorphe, on peut supposer que 
$$
P=\sum_{k=1}^pX_k(x)\wedge \frac{\partial}{\partial x_{n+k}} + \sum_{1\leq i<j\leq n} \{x_i,x_{j}\}\frac{\partial}{\partial x_i}\wedge \frac{\partial}{\partial x_{j}}.
$$
v\'erifie les propri\'et\'es suivantes, $U$ \'etant un voisinage born\'e de l'origine
dans $\Bbb C^p$~:
\begin{itemize}
\item pour $1\leq k\leq p$, $X_k$ est un champ de vecteurs holomorphe de
  $\Bbb C^n$ d\'ependant holomorphiquement des $p$ param\`etres $x''$,
  i.e. $X_k\in {\cal O}_p(\bar U)\otimes \vfo n 1$;
\item les $X_k$ commutent deux \`a deux, i.e. $[X_i,X_j]=0$ pour $1\leq
  i,j\leq p$;
\item pour $1\leq k\leq p$, $X_k\equiv S_k \mod{{\cal O}_p(\bar U)\otimes \vfo n 2}$;
\item pour $1\leq i,j \leq n$, $\{x_i,x_{j}\}\in {\cal O}_p(\bar U)\otimes {\cal M}^2$.
\end{itemize}
La partie lin\'eaire de $P$ \`a l'origine est alors 
$$
{\cal L}:= \sum_{k=1}^pS_k(x')\wedge \frac{\partial}{\partial x_{n+k}}.
$$
On rappelle que l'on a pos\'e, pour $1\leq i\leq n$, 
$$
\Lambda_i=\sum_{j=1}^p\lambda_{j,i}\frac{\partial}{\partial x_{n+j}}.
$$
Puisque la famille $S$ est de rang $p$, {\bf on peut supposer que la famille
$\{\Lambda_i\}_{1\leq i\leq p}$ est libre sur $\Bbb C$, ce que nous ferons
dans la suite}.

pour $k\in \Bbb N^*$, on pose
$$
\omega_k(S)=\inf \left\{\max_i|(Q,\lambda^i)-\lambda_{i,j}|\neq 0, Q\in \Bbb N^n, 2\leq |Q|\leq 2^k, 1\leq j\leq n\right\}
$$
Nous dirons que la famille $S$ est {\bf diophantienne} si la condition arithm\'etique
$$
\omega(S)~:\quad -\sum_{k\geq 1}\frac{\log\omega_k(S) }{2^k}<+\infty.
$$
est satisfaite.
\begin{rem}
Lorsque la famille $S$ ne contient qu'un seul \'el\'ement, la condition $\omega(S)$ n'est autre que la condition de Bruno \cite{bruno} (qui est plus faible que celle de Siegel\cite{Arn2}[p. 182]).
\end{rem}
\begin{defi}
Un ensemble compact $K$ sera dit {\bf bon} s'il est de Stein et si, pour tout
ensemble analytique $V$ d'un voisinage de $K$, $V\cap K$ n'a qu'un nombre fini
de composantes connexes. Un ensemble compact $K$ sera dit compact de $0$ s'il
contient un voisinage de $0$. Un polydisque ferm\'e est un bon compact de son centre.
\end{defi}
\begin{rem}
Dans la pratique, nous pourrions ne consid\'erer, comme bon compact, que le cas d'un point ou d'un
polydisque ferm\'e centr\'e en ce point. La propri\'et\'e essentielle que nous
utiliserons est la noetherianit\'e de l'anneau des fonctions holomorphes au
voisinage d'un bon compact \cite{frisch, siu-noether}. Nous renvoyons le lecteur aux
ouvrages \cite{hormanderscv, malgrangescv} pour les notions sur les espaces de Stein.
\end{rem}
\begin{theo}\label{theo1}
Supposons que $K:=\bar U$ soit un bon compact de $0\in \Bbb C^p$.
Supposons que la famille $S=\{S_i\}_{1\leq i\leq
  p}$ des parties lin\'eaires soit {\bf
  diophantienne}. Supposons que la famille commutative $\{X_i\}_{1\leq i\leq
  p}$ admette un forme normale formelle de la forme $\{\sum_{j=1}^p \hat
  a_{i,j}S_j\}_{1\leq i\leq p}$ o\`u les $\hat a_{i,j}$ appartiennent \`a
  l'anneau ${\cal O}_p(K)\otimes\widehat{\cal O}_n^S$. Alors, il existe un
  diff\'eomorphisme holomorphe $\Phi$ de $(\Bbb C^N,0)$ dans lui-m\^eme, fixant
  l'origine et tangent \`a l'identit\'e en ce point tel que 
\begin{eqnarray*}
\Phi_*P & = & \sum_{k=1}^p\left(\sum_{l=1}^p\tilde a_{k,l}(x)S_l\right)\wedge \frac{\partial }{\partial x_{n+k}}+\sum_{p+1\leq i<j\leq n}c_{i,j}x_ix_j\frac{\partial }{\partial x_{i}}\wedge \frac{\partial }{\partial x_{j}}\\
& & +\sum_{1\leq i<j\leq n}\left(\sum_{\substack{Q\in \Bbb N^n_2,\\
Q\neq E_i+E_j\\ \forall r\;
(Q,\lambda^r)=\lambda_{r,i}+\lambda_{r,j}}}g_{i,j,Q}(x'')(x')^{Q}\right)\frac{\partial
}{\partial x_{i}}\wedge \frac{\partial }{\partial x_{j}}
\end{eqnarray*}
o\`u les $\tilde a_{r,s}$ appartiennent \`a ${\cal O}_p(K')\otimes {\cal
  O}_n^S$, les $g_{i,j,Q}$ appartiennent \`a $ {\cal O}_p(V')$, $K'$
  (resp. $V'$) \'etant
  un bon compact de $0$ dans $\Bbb C^p$ (resp. voisinage de $K'$), et les $c_{i,j}$ sont des nombres
  complexes. 
\end{theo}

Les corollaires suivants sont des applications imm\'ediates du th\'eor\`eme.
\begin{defi}
Nous dirons que la partie lin\'eaire ${\cal L}$ de $P$ est {\bf non-r\'esonnante} si, pour tout
$n$-uplet d'entiers $(q_1,\ldots, q_n)$ tel que $q_i\geq -1$, deux au plus
pouvant \^etre simultan\'ement \'egaux \`a $-1$, il existe un entier $1\leq
j\leq p$ tel que 
$$
q_1\lambda_{j,1}+\ldots+q_n\lambda_{j,n}\neq 0.
$$
On appellera {\bf relation de r\'esonnance} de ${\cal L}$, une \'egalit\'e de la forme
$$
q_1\lambda_{j,1}+\ldots+q_n\lambda_{j,n}= \lambda_{j,r}+\lambda_{j,s}\quad\quad \forall 1\leq j\leq n,
$$
o\`u $Q=(q_1,\ldots, q_n)\in \Bbb N^n$.
La famille $S$ est non-r\'esonnante lorsque $(\fvfo n 2)^S=\{0\}$.
\end{defi}
\begin{rem}
Si ${\cal L}$ est non-r\'esonnante alors la famille $S$ est non-r\'esonnante.
\end{rem}
\begin{coro}
Supposons que la famille $S=\{S_i\}_{1\leq i\leq
  p}$ des parties lin\'eaires soit {\bf
  diophantienne} et que ${\cal L}$ soit {\bf non-r\'esonnante}. 
Alors alors il existe un diff\'eomorphisme holomorphe $\Phi$ de $(\Bbb C^N,0)$ dans lui-m\^eme, fixant l'orgine et tangent \`a l'identit\'e en ce point tel que 
\begin{eqnarray*}
\Phi_*P & = & \sum_{r=1}^p S_r\wedge \frac{\partial }{\partial x_{n+r}}+\sum_{p+1\leq i<j\leq n}c_{i,j}x_ix_j\frac{\partial }{\partial x_{i}}\wedge \frac{\partial }{\partial x_{j}};
\end{eqnarray*}
o\`u les  $c_{i,j}$ sont des constantes.
\end{coro}

Ce dernier corollaire, dans le cas $p=1$, est du \`a J.-P. Dufour et
M. Zhitomirskii \cite{dufour3}[theorem 6.1]. En particulier, si $n\leq p+1$, la structure de Poisson est lin\'earisable (sans condition de rang).

\begin{theo}\label{theo2}
Sous les hypoth\`eses du th\'eor\`eme pr\'ec\'edent, supposons de plus, que le
rang de $P$ soit egal \`a $2p$ avec $n>p+1$. Alors il existe un diff\'eomorphisme holomorphe $\Phi$ de $(\Bbb C^N,0)$ dans lui-m\^eme, fixant l'orgine et tangent \`a l'identit\'e en ce point tel que 
$$
\Phi_*P  =  \sum_{k=1}^p\left(\sum_{l=1}^p b_{k,l}(x^{R_1},\ldots, x^{R_t})S_l\right)\wedge \frac{\partial }{\partial x_{n+k}}.
$$
o\`u les $b_{k,l}$ appartiennent \`a ${\cal O}_n^S$.
\end{theo}

\begin{rem}
Le r\'esultat pr\'ec\'edent admet {\bf l'interpr\'etation g\'eom\'etrique} suivante. On supposera que l'anneau des invariants $\widehat{\cal O}_n^S$ n'est pas r\'eduit aux constantes. Soit $\pi: (\Bbb C^n,0)\rightarrow (\Bbb C^t,0)$ l'application d\'efinie par
$\pi(x)= (x^{R_1},\ldots,x^{R_t})$ avec $\widehat{\cal O}_n^S=\Bbb
C[[x^{R_1},\ldots,x^{R_t}]]$ (voir section \ref{forme-normale}). Il existe
alors un voisinage ${\cal U}$ de l'origine dans $\Bbb C^n$ tel que, pour tout
$b\in \pi({\cal U})$, la structure de Poisson $P$ induit, sur la vari\'et\'e
torique "g\'en\'eralis\'ee" $(\pi^{-1}(b)\cap {\cal U})\times \Bbb C^p$, la
structure de Poisson lin\'eaire 
$$
\sum_{r=1}^p\left(\sum_{s=1}^p b_{r,s}(b)S_s\right)\wedge \frac{\partial
  }{\partial x_{n+r}}\left|_{(\pi^{-1}(b)\cap {\cal U})\times \Bbb C^p}\right..
$$ 
En effet, par d\'efinition, les champs $S_k$ sont tangents aux fibres $(\pi^{-1}(b)\cap {\cal U})$.
\end{rem}
\begin{coro}
Soit $P$ une structure de Poisson v\'erifiant les hypoth\`eses $(H)$ et de
rang egal \`a $2p$, $n>p+1$. On suppose que la famille des hamiltoniens associ\'es aux variables $x''$ est {\bf formellement lin\'earisable} et que la famille de leurs parties lin\'eaires \`a l'origine est {\bf diophantienne}. Alors la structure de Poisson $P$ est holomorphiquement lin\'earisable.
\end{coro}

Ce r\'esultat est du \`a J.-P. Dufour \cite{dufour2}, dans le cadre ${\cal C}^{\infty}$, lorsque la famille $S$ est {\bf hyperbolique et
  non-r\'esonnante}. 
\begin{rem}
Tous ces r\'esultats sont vrais dans la cat\'egorie analytique r\'eelle.
\end{rem}


\section{Forme normale de famille commutative de champs de vecteurs}

L'objet de cette section est de rappeler un des r\'esultats de l'auteur \cite{stolo-ihes} que nous utiliserons dans cet article.
Soit $\{X_1,\ldots, X_p\}$ une famille commutative de champs de vecteurs holomorphes au voisinage de l'origine dans $\Bbb C^n$. On les suppose nuls en ce point et on suppose que leur partie lin\'eaire sont diagonales et ind\'ependantes sur $\Bbb C$. On posera $J^1(X_i)=S_i:=\sum_{j=1}^n\lambda_{i,j}x_j{\partial }/{\partial x_j}$, $\lambda^i=(\lambda_{i,1},\ldots, \lambda_{i,n})$ et $S$ d\'esignera la famille des $S_i$.

On a alors le r\'esultat suivant ~:
\begin{theo}\cite{stolo-ihes}
Supposons que la famille des parties lin\'eaires $\{S_i\}_{1\leq i\leq p}$ soit {\bf diophantienne}. Supposons que la famille commutative $\{X_i\}_{1\leq i\leq p}$ admette un forme normale formelle de la forme $\{\sum_{j=1}^p \hat a_{i,j}S_j\}_{1\leq i\leq p}$ o\`u les $\hat a_{i,j}$ appartiennent \`a l'anneau $\widehat{\cal O}_n^S$.
Alors il existe un germe de diff\'eomorphisme holomorphe $\Phi$ de $(\Bbb C^n,0)$ dans lui-m\^eme, tangent \`a l'identit\'e en $0$, tel que, pour $1\leq i\leq p$
$$
\Phi_*X_i=\sum_{j=1}^p a_{i,j}S_j,
$$
o\`u les $a_{i,j}\in \widehat{\cal O}_n^S\cap{\cal O}_n$.
\end{theo}

Dans ce qui suit, nous aurons besoin de la "version \`a param\`etres" suivante :
\begin{theo}\label{param}
Soit $K$ un bon compact de $0\in\Bbb C^k$. On notera $x''$ les coordonn\'ees de $\Bbb C^k$. Soit $\{X_i\}_{1\leq i\leq p}$ une famille commutative d'\'el\'ements de ${\cal O}_k(K)\otimes \vfo n 1$. Soit $S_i$ la partie lin\'eaire de $X_i$ en $0\in \Bbb C^n$, que l'on suppose ind\'ependante de $x''$. Supposons que la famille $\{S_i\}_{1\leq i\leq p}$ soit {\bf diophantienne}.
Supposons, en outre, que la famille commutative $\{X_i\}_{1\leq i\leq p}$ admette une forme normale formelle de la forme $\{\sum_{j=1}^p \hat a_{i,j}S_j\}_{1\leq i\leq p}$ o\`u les $\hat a_{i,j}$ appartiennent \`a l'anneau ${\cal O}_k(K)\otimes\widehat{\cal O}_n^S$.
Alors il existe un germe de diff\'eomorphisme holomorphe $\Phi$ de $(\Bbb C^n,0)$ dans lui-m\^eme, tangent \`a l'identit\'e en $0$, d\'ependant holomorphiquement de $x''$ dans $K$ et tel que , pour $1\leq i\leq p$
$$
\Phi_*X_i=\sum_{j=1}^p a_{i,j}S_j,
$$
o\`u les $a_{i,j}\in {\cal O}_k(K)\otimes\left(\widehat{\cal O}_n^S\cap{\cal O}_n\right)$.
\end{theo}
En effet, supposons que $\hat \Phi$ conjugue $X_j$ \`a $Y_j$, pour $1\leq j\leq p$, avec
\begin{eqnarray*}
\hat \Phi(x, x'')& = & \left(x_i+\sum_{Q\in \Bbb N^n, |Q|\geq 2}\phi_{i,Q}(x'')x^Q\right)_{1\leq i\leq n},\\
X_j & = & S_j +\sum_{i=1}^n\left(\sum_{Q\in \Bbb N^n, |Q|\geq 2}X_{j,i,Q}(x'')x^Q\right)\frac{\partial}{\partial x_i},\\
\text{et }\;Y_j & = & S_j +\sum_{i=1}^n\left(\sum_{Q\in \Bbb N^n, |Q|\geq 2}Y_{j,i,Q}(x'')x^Q\right)\frac{\partial}{\partial x_i}.
\end{eqnarray*}

Les \'equations que doivent satisfaire ces coefficients sont de la forme
$$
\left((Q,\mu^j)-\mu_{j,i}\right)\phi_{i,Q}(x'')+Y_{j,i,Q}(x'')= X_{j,i,Q}(x'')+ P_{j,i,Q}(x''),
$$
o\`u $P_{j,i,Q}$ est un polyn\^ome en les $\phi_{i,Q'}(x'')$, $Y_{j,i,Q'}(x'')$, $X_{j,i,Q'}(x'')$, $|Q'|<|Q|$. La d\'efinition de $\phi_{i,Q}(x'')$ et de $Y_{j,i,Q}(x'')$ se fait par r\'eccurence comme suit : si pour un $j$, on a $(Q,\mu^j)\neq\mu_{j,i}$ alors on pose $Y_{j,i,Q}(x'')=0$ et $\phi_{i,Q}(x'')$ est d\'etermin\'e par la $j$-i\`eme \'equation. Sinon, on pose 
$\phi_{i,Q}(x'')=0$ et $Y_{j,i,Q}(x'')$ est d\'etermin\'e par le second membre. 
Par cons\'equent, on montre par r\'ecurrence que si les $X_{j,i,Q}(x'')$ sont
holomorphes sur un m\^eme voisinage $V$ de $K$, il en est de m\^eme des $\phi_{i,Q}(x'')$ et des
$Y_{j,i,Q}(x'')$.


\section{Preuve du th\'eor\`eme \ref{theo1}}

Utilisons le th\'eor\`eme \ref{param}, pour normaliser la famille commutative 
de champs de vecteurs $\{X_i\}_{1\leq i\leq p}$. Dans ces nouvelles
coordonn\'ees, d\'eveloppons l'identit\'e de Jacobi, pour $1\leq i,j\leq n$ et
$1\leq m\leq p$, 
\begin{eqnarray}\label{jacobi}
\{x_{n+m},\{x_i,x_j\}\} & = & -\{x_{i},\{x_j,x_{n+m}\}\}+\{x_{j},\{x_i,x_{n+m}\}\}.
\end{eqnarray}
Tout d'abord, on a 
\begin{eqnarray*}
\{x_{n+m},\{x_i,x_j\}\} & = & \sum_{k=1}^N\{x_k, x_{n+m}\}\frac{\partial
  \{x_i,x_j\}}{\partial x_k}\\
& = & \sum_{k=1}^n\{x_k, x_{n+m}\}\frac{\partial \{x_i,x_j\}}{\partial x_k}=
  X_m(\{x_i,x_j\}).
\end{eqnarray*}
On a donc
\begin{eqnarray}\label{terme1}
\{x_{n+m},\{x_i,x_j\}\} & = & \sum_{k=1}^p a_{m,k}(x)S_k(\{x_i,x_j\}).
\end{eqnarray}
D'autre part, par d\'efinition, on a $\{x_j,x_{n+m}\}=X_{m,j}$. On a alors 
\begin{eqnarray*}
\{x_{i},\{x_j,x_{n+m}\}\} & = & \sum_{k=1}^N\{x_k, x_{i}\}\frac{\partial
  X_{m,j}}{\partial x_k}\\
& = & \sum_{k=1}^n\{x_k, x_{i}\}\frac{\partial X_{m,j}}{\partial x_k} +
  \sum_{k=1}^p\{x_{n+k}, x_{i}\}\frac{\partial X_{m,j}}{\partial x_{n+k}}\\
& = & \sum_{k=1}^n\{x_k, x_{i}\}\frac{\partial X_{m,j}}{\partial x_k} -
  \sum_{k=1}^p X_{k,i}\frac{\partial X_{m,j}}{\partial x_{n+k}}.
\end{eqnarray*}
Or
\begin{eqnarray}\label{terme2}
\sum_{k=1}^n\{x_k, x_{i}\}\frac{\partial X_{m,j}}{\partial x_k} & = & 
\{x_i,x_j\}\left(\sum_{k=1}^p a_{m,r} \lambda_{r,j}\right)\\
& & +x_j\sum_{k=1}^n \{x_k,x_i\}\frac{\partial \sum_{r=1}^p a_{m,r}
  \lambda_{r,j}}{\partial x_k}
\end{eqnarray}
Consid\'erons alors l'op\'erateur ${\cal }O_p(K)\otimes {\cal O}_n$-lin\'eaire 
$$
N : {\cal  M}_n\left({\cal O}_p(K)\otimes {\cal O}_n\right) \rightarrow {\cal
  M}_n\left({\cal O}_p(K)\otimes {\cal O}_n\right)\otimes \left({\cal O}_p(K)\otimes
  {\cal O}_n\right)^n
$$ 
d\'efini par 
$$
N\left((f_{i,j})_{1\leq i,j\leq n}\right)= (N_{i,j}(f))_{1\leq i,j\leq n}=
\left(
  \left(\sum_{r=1}^p\sum_{k=1}^n\left(x_if_{k,j}\lambda_{r,i}-x_jf_{k,i}\lambda_{r,j}\right)\frac{\partial a_{m,r}}{\partial x_k}\right)_{1\leq m\leq n}\right)_{1\leq i,j\leq n}.
$$
On posera
\begin{eqnarray}\label{nij}
N_{i,j,m}(f) :=
\sum_{r=1}^p\sum_{k=1}^n\left(x_if_{k,j}\lambda_{r,i}-x_jf_{k,i}\lambda_{r,j}\right)\frac{\partial
  a_{m,r}}{\partial x_k}.
\end{eqnarray}
Soit 
\begin{eqnarray*}
M_{i,j,m} & := & \sum_{l=1}^p X_{l,i}\frac{\partial X_{m,j}}{\partial
  x_{n+l}}-\sum_{l=1}^p X_{l,j}\frac{\partial X_{m,i}}{\partial x_{n+l}}\\
& = & x_ix_j\sum_{l=1}^p\left[\left(\sum_{k=1}^p a_{l,k}
  \lambda_{k,i}\right)\left(\frac{\partial \sum_{k=1}^p a_{m,k}
  \lambda_{k,j}}{\partial x_{l+m}}\right)-\left(\sum_{k=1}^p a_{l,k}
  \lambda_{k,j}\right)\left(\frac{\partial \sum_{k=1}^p a_{m,k}
  \lambda_{k,i}}{\partial x_{n+l}}\right)\right].
\end{eqnarray*}

\begin{lemm}
Pour $1\leq i,j\leq n$ et $1\leq m\leq p$, 
\begin{itemize}
\item $M_{i,j,m}$ appartient \`a
l'id\'eal engendr\'e par $x_ix_j$ dans ${\cal O}_p(K)\otimes {\cal O}_n^S$,
\item sous les hypoth\`eses $(H)$, $M_{i,j,m}$ est d'ordre sup\'erieur ou
  \'egal \`a $5$,\item si chaque $f_{i,j}$ est de la forme 
$$
\sum_{\substack{Q\in \Bbb N^n_2,\\ \forall
    r,\;(Q,\lambda^r)=\lambda_{r,i}+\lambda_{r,j}}}f_{i,j,Q}(x'')(x')^Q
$$
avec $f_{i,j,Q}\in {\cal O}_p(V)$, $V$ \'etant un voisinage de $K$; il en est de m\^eme pour $N_{i,j,m}(f)$.
\end{itemize}
\end{lemm}
\begin{proof}
Le premier point est clair car les fonctions $a_{i,j}$ appartiennent \`a
${\cal O}_p(K)\otimes {\cal O}_n^S$; leurs d\'eriv\'ees par rapport aux
variables $x''$ aussi. Les fonctions $a_{i,j}$ sont de la forme
$$
\sum_{\substack{P\in \Bbb N^n_2,\\\forall r\; (P,\lambda^r)=0}}a_{i,j,P}(x'')(x')^P.
$$
D'apr\`es les hypoth\`eses $(H)$, les mon\^omes "r\'esonnants" qui engendrent
${\cal O}_n^S$ sont de degr\'e sup\'erieur ou \'egal \`a $3$.
Par cons\'equent, la fonction $x_if_{k,j}{\partial a_{m,r}}/{\partial x_k}$
est de la forme 
$$
\sum_{\substack{Q\in \Bbb N^n_2,\\ \forall r\; (Q,\lambda^r)=\lambda_{r,k}+\lambda_{r,j},\\
    (P,\lambda^r)=0}}g_{i,j,Q}(x'')(x')^{Q+E_i+P-E_k},
$$
c'est-\`a-dire de la forme 
$$
\sum_{\substack{R\in \Bbb N^n_3,\\ \forall r\;
    (R,\lambda^r)=\lambda_{r,i}+\lambda_{r,j}}}\tilde g_{i,j,R}(x'')(x')^{R}.
$$
En fait, on a m\^eme, 
$$
x_if_{k,j}\frac{\partial a_{m,r}}{\partial x_k}=x_i\sum_{\substack{R\in \Bbb N^n_2,\\ \forall r\;
    (T,\lambda^r)=\lambda_{r,j}}}\tilde g_{i,j,T}(x'')(x')^{T}.
$$
\end{proof}

L'\'equation de Jacobi $(\ref{jacobi})$ s'\'ecrit, pour $1\leq m \leq p$,
$$
X_m(\{x_i,x_j\})-\{x_i,x_j\}\left(\sum_{r=1}^p a_{m,r}(\lambda_{r,j}+\lambda_{r,i})\right)=M_{i,j,m}+N_{i,j,m}\left((\{x_r,x_s\})_{1\leq r,s\leq n}\right).
$$
Soit $A=(a_{k,l})_{1\leq k,l\leq p}$; par hypoth\`ese, $A(0, x'')=Id$ et $a_{k,l}\in {\cal O}_p(K)\otimes {\cal O}_n^S$. La matrice $A$ est donc
inversible au voisinage de $\{0\}\times K$ dans $\Bbb C^n\times\Bbb C^p$ et les coefficients de $A^{-1}$
appartiennent \`a  ${\cal O}_p(K)\otimes {\cal O}_n^S$. Puisque l'on a
$X_m=\sum_{k=1}^pa_{m,k}S_k$, les $p$ \'equations pr\'ec\'edentes peuvent
s'\'ecrire sous la forme matricielle 
\begin{eqnarray}\label{jacobi2}
\begin{pmatrix} S_1(\{x_i,x_j\}) \\ \vdots \\ S_p(\{x_i,x_j\})\end{pmatrix} -
\{x_i,x_j\}\begin{pmatrix} \lambda_{1,i}+\lambda_{1,j} \\ \vdots \\
  \lambda_{p,i}+\lambda_{p,j}\end{pmatrix} & = &
  \tilde N_{i,j}\left((\{x_r,x_s\})_{1\leq r,s\leq n}\right)+ \tilde M_{i,j}\\
\end{eqnarray}
o\`u l'on a pos\'e  $\tilde N_{i,j}= A^{-1}N_{i,j}$ et $\tilde
M_{i,j}=A^{-1}M_{i,j}$. 
Les composantes de $\tilde M_{i,j}$ appartiennent \`a l'id\'eal engendr\'e par
$x_ix_j$ dans ${\cal O}_p(K)\otimes {\cal O}_n^S$. Si chancun des $f_{i,j}$
v\'erifie les hypoth\`eses du lemme pr\'ec\'edent, la conclusion reste valide
pour les $\tilde N_{i,j}(f)$.
\begin{lemm}
Pour $1\leq i,j\leq n$, on a 
$$
\{x_i,x_j\}= \sum_{\substack{Q\in \Bbb N^n_2,\\ \forall r\;
    (Q,\lambda^r)=\lambda_{r,i}+\lambda_{r,j}}}g_{i,j,Q}(x'')(x')^{Q},
$$ 
o\`u les $g_{i,j,Q}$ appartiennent \`a ${\cal O}_p(V)$, $V$ \'etant un voisinage de $K$. De plus, si l'on
suppose que, pour $1\leq i,j\leq n$, les solutions $Q\in \Bbb N^2$,
$2\leq |Q|$, des \'equations $(Q,\lambda^r)=\lambda_{r,i}+\lambda_{r,j}$ pour $1\leq r\leq p$ sont de la forme $Q=R+E_i+E_j$, $R\in \Bbb N^n_2$ alors
$$
\{x_i,x_j\}=x_ix_j \sum_{\substack{Q\in \Bbb N^n_2,\\ \forall r\;
    (Q,\lambda^r)=0}}g_{i,j,Q}(x'')(x')^{Q};
$$ 
c'est-\`a-dire que $\{x_i,x_j\}$ appartient \`a l'id\'eal de ${\cal  O}_p(K)\otimes {\cal O}_n^S$ engendr\'e par $x_ix_j$.
\end{lemm}
\begin{proof}
Posons 
$$
\{x_i,x_j\}= \sum_{Q\in \Bbb N^n_2} g_{i,j,Q}(x'')(x')^Q.
$$
On rappelle, d'apr\`es le lemme \ref{ordre2}, que $\{x_i,x_j\}$ est d'ordre
sup\'erieur ou \'egal \`a deux en $(0,x'')$. On a alors, pour $1\leq k\leq p$ 
$$
S_k(\{x_i,x_j\})-\{x_i,x_j\}(\lambda_{k,i}+\lambda_{k,j}) = \sum_{Q\in \Bbb N^n,
  |Q|\geq 2} \left((Q,\lambda^k)-(\lambda_{k,i}+\lambda_{k,j})\right)g_{i,j,Q}(x'')(x')^Q.
$$
D'autre part, on a, pour tout entier $d\geq 1$
$$
J^d\left(\tilde N_{i,j,m}(f)\right)= J^d\left(\tilde
  N_{i,j,m}\left(J^{d-1}(f)\right)\right),
$$
o\`u $J^{d-1}(f)$ d\'esigne la matrice dont les composantes sont les
$J^{d-1}(f_{i,j})$ (les jets sont bien entendu en $(0,x'')$ en les variables
$x'$). 
En utilisant la d\'efinition $(\ref{nij})$ de $N_{i,j,m}$, il est clair que si
$f$ est d'ordre sup\'erieur ou \'egal \`a $d$ alors $N_{i,j,m}(f)$ est d'ordre
sup\'erieur ou \'egal \`a $d+1$ (en fait, sous les hypoth\`eses $(H)$, l'ordre
de $N_{i,j,m}(f)$ est sup\'erieur ou \'egal \`a $d+5$). La multiplication par la matrice $A^{-1}$ ne change rien quant \`a l'ordre.
D'autre part, pour les m\^emes raisons, $M_{i,j,m}$ est d'ordre sup\'erieur ou \'egal \`a cinq. Il en est de m\^eme pour $\tilde M_{i,j,m}$.

Montrons alors le r\'esultat par r\'ecurrence sur la longueur du multiindice
$|Q|\geq 2$. Pour $|Q|=2$, l'analyse faite pr\'ec\'edement montre que la composante
de $(x')^Q$ du second membre de l'identit\'e de Jacobi $(\ref{jacobi2})$ est
nulle. Lorsqu'il existe $1\leq r\leq p$ tel que
$(Q,\lambda^r)\neq\lambda_{r,i}+\lambda_{r,j}$ alors $g_{i,j,Q}\equiv 0$. L'assertion est
donc vraie.
Supposons la vraie jusqu'\`a l'ordre $d$. D'apr\`es ce qui pr\'ec\`ede, le coefficient de $(x')^Q$, $|Q|=d+1$, du second membre de l'identit\'e de Jacobi $(\ref{jacobi2})$ est la somme de celui de $\tilde M_{i,j}$ et de celui de $\tilde N_{i,j,m}\left(J^{d-1}\left(\left(\{x_r,x_s\}\right)_{1\leq r,s\leq n}\right)\right)$. D'apr\`es le lemme pr\'ec\'edent, s'il existe $1\leq r\leq p$ tel que $(Q,\lambda^r)\neq\lambda_{r,i}+\lambda_{r,j}$, alors ce coefficient est nul. On a alors $\left((Q,\lambda^r)-\lambda_{r,i}-\lambda_{r,j}\right)g_{i,j,Q}\equiv 0$; c'est-\`a-dire $g_{i,j,Q}\equiv 0$. Ce qui ach\`eve la r\'ecurrence.

En ce qui concerne le deuxi\`eme point, il suffit de remarquer que si toute solution de $(Q,\lambda^r)=\lambda_{r,i}+\lambda_{r,j}$, pour $1\leq r\leq p$, est de la forme $Q=R+E_i+E_j$ alors $(R,\lambda^r)=0$ pour $1\leq r\leq p$.
\end{proof}

Nous avons donc montr\'e que, dans un bon syst\`eme de coordonn\'ees holomorphes, on peut \'ecrire
$$
P=\sum_{r=1}^p\left(\sum_{s=1}^pa_{r,s}(x)S_s\right)\wedge \frac{\partial }{\partial x_{n+r}}+\sum_{1\leq i<j\leq n}\left(\sum_{\substack{Q\in \Bbb N^n_2,\\ \forall r\; (Q,\lambda^r)=\lambda_{r,i}+\lambda_{r,j}}}g_{i,j,Q}(x'')(x')^{Q}\right)\frac{\partial }{\partial x_{i}}\wedge \frac{\partial }{\partial x_{j}},
$$
o\`u $a_{r,s}\in {\cal O}_p(K)\otimes {\cal O}_n^S$ et $g_{i,j,Q}\in {\cal
  O}_p(V)$, $V$ \'etant un voisinage de $K$.

Soient $\beta_1(x''),\ldots, \beta_n(x'')$ des fonctions holomorphes au voisinage de $K$ telles que $\beta_i(0)=1$.

Posons alors $x_i=\beta_i(y'')y_i$, $i=1,\ldots, n$, et $x_{n+k}=y_{n+k}$ pour
$k=1,\ldots, p$. On notera $x=\Psi(y)$ cette transformation tangente \`a
l'identit\'e en $0$. Dans ces nouvelles coordonn\'ees, le crochet de Poisson $\{.,.\}'$ est d\'efini par 
$$
\{f',g'\}'(y)=\{f'\circ \Psi^{-1}, g'\circ \Psi^{-1}\}\circ \Psi(y).
$$
En particulier, on a, pour $1\leq i,j\leq n$
\begin{eqnarray*}
\{y_j, y_i\}' & = & \left\{\frac{x_j}{\beta_j(x'')},
  \frac{x_i}{\beta_i(x'')}\right\}\circ \Psi(y)\\
& = &\left(\frac{1}{\beta_j(x'')}\left\{x_j,
  \frac{x_i}{\beta_i(x'')}\right\}+x_j\left\{\frac{1}{\beta_j(x'')},
  \frac{x_i}{\beta_i(x'')}\right\}\right)\circ \Psi(y)\\
& = & \left(\frac{1}{\beta_j(x'')}\left(\frac{1}{\beta_i(x'')}\{x_j,
  x_i\}+x_i\left\{x_j,\frac{1}{\beta_i(x'')}\right\}\right)\right.\\
& &
  \left.+\frac{x_j}{\beta_i(x'')}\left\{\frac{1}{\beta_j(x'')},x_i\right\}+x_ix_j\left\{\frac{1}{\beta_j(x'')},\frac{1}{\beta_i(x'')}\right\}\right)\circ
  \Psi(y)\\
& = & \left(\frac{\{x_j,x_i\}}{\beta_i(x'')\beta_j(x'')}+\frac{x_i}{\beta_j(x'')}\left\{x_j,\frac{1}{\beta_i(x'')}\right\}-\frac{x_j}{\beta_i(x'')}\left\{x_i,\frac{1}{\beta_j(x'')}\right\}\right)\circ \Psi(y)\\
\end{eqnarray*}
gr\^ace \`a l'identit\'e de Leibnitz et au fait que $\{x_{n+k},x_{n+l}\}=0$ pour $1\leq k,l\leq p$. Par d\'efinition, on a 
\begin{eqnarray*}
\left\{x_i,\frac{1}{\beta_j(x'')}\right\}& = &
\frac{-1}{\beta_i^2(x'')}\sum_{k=1}^p\{x_{n+k}, x_i\}\frac{\partial
  \beta_j}{\partial x_{n+k}}\\
& = & \frac{1}{\beta_j^2(x'')}\sum_{k=1}^p X_{k,i}\frac{\partial
  \beta_j}{\partial x_{n+k}}\\
& = & x_i\frac{\Lambda_i(\beta_j)}{\beta_j^2(x'')}\quad\mod{{\cal O}_p(K)\otimes
  {\cal M}_n^2}.\\
\end{eqnarray*}
Par cons\'equent, le coefficient de $y_iy_j$ dans $\{y_j, y_i\}'$ est \'egal
\`a 
$$
g_{j,i,E_i+E_j}(x'')+\frac{\Lambda_j(\beta_i)}{\beta_i(x'')}-\frac{\Lambda_i(\beta_j)}{\beta_j(x'')}.
$$
En posant, $\beta_i(x'')=\exp \gamma_i$, ce coefficient s'\'ecrit
\begin{eqnarray}\label{cobord}
g_{j,i,E_i+E_j}(x'')+\Lambda_j(\gamma_i)-\Lambda_i(\gamma_j).
\end{eqnarray}

\begin{lemm}
Pour $1\leq i,j,k\leq n$, on a 
$$
\Lambda_j(g_{k,i,E_k+E_i})-\Lambda_i(g_{k,j,E_k+E_j})+\Lambda_k(g_{j,i,E_j+E_i})=0.
$$
\end{lemm}
\begin{proof}
Consid\'erons l'identit\'e de Jacobi
$$
\{x_k,\{x_i,x_j\}\}+\{x_i,\{x_j,x_k\}\}+\{x_j,\{x_k,x_i\}\}=0.
$$
On a 
\begin{eqnarray*}
\{x_k,\{x_i,x_j\}\} & = & \sum_{m=1}^N\{x_m,x_k\}\frac{\partial \{x_i,x_j\}}{\partial x_m}\\
& = & \sum_{m=1}^n\{x_m,x_k\}\frac{\partial \{x_i,x_j\}}{\partial x_m}-\sum_{l=1}^pX_{l,k}\frac{\partial \{x_i,x_j\}}{\partial x_{n+l}}. 
\end{eqnarray*}
On a, pour $1\leq m\leq n$ 
$$
\{x_m,x_k\}\frac{\partial \{x_i,x_j\}}{\partial x_m}=\left(\sum_{\substack{P\in \Bbb N^n_2\\  (P,\lambda)=\lambda_{k}+\lambda_{m}}}g_{m,k,P}(x')^P\right)\left(\sum_{\substack{Q\in \Bbb N^n_2,\\  (Q,\lambda)=\lambda_{i}+\lambda_{j}}}q_mg_{i,j,Q}(x')^{Q-E_m}\right)
$$
o\`u, afin d'all\'eger l'\'ecriture, $(Q,\lambda)=\lambda_{i}+\lambda_{j}$
signifie dans les sommes: {\it pour tout $1\leq t\leq p$,
  $(Q,\lambda^t)=\lambda_{t,i}+\lambda_{t,j}$}. On a aussi, pour $1\leq l\leq p$
$$
X_{l,k}\frac{\partial \{x_i,x_j\}}{\partial x_{n+l}} = x_k\left(\sum_{s=1}^pa_{l,s}\lambda_{s,k}\right)\left(\sum_{\substack{P\in \Bbb N^n_2\\  (P,\lambda)=\lambda_{i}+\lambda_{j}}}\frac{\partial g_{i,j,P}}{\partial x_{l+n}}(x')^P\right).
$$

Prenons le jet d'ordre $3$ en $0\in \Bbb C^n$(i.e. en x') de ces expressions. Il vient  
$$
J^3\left(\{x_m,x_k\}\frac{\partial \{x_i,x_j\}}{\partial x_m}\right)=\left(\sum_{\substack{ \lambda_{r}+\lambda_{s}=\lambda_{k}+\lambda_{m}}}g_{m,k,E_r+E_s}x_rx_s\right)\left(\sum_{\substack{\lambda_{r}+\lambda_{m}=\lambda_{i}+\lambda_{j}}}g_{i,j,E_r+E_m}x_r\right)
$$
et
$$
J^3\left(X_{l,k}\frac{\partial \{x_i,x_j\}}{\partial x_{n+l}}\right) = \lambda_{l,k}x_k\left(\sum_{\substack{\lambda_{r}+\lambda_{s}=\lambda_{i}+\lambda_{j}}}\frac{\partial g_{i,j,E_r+E_s}}{\partial x_{l+n}}x_rx_s\right).
$$
Le coefficient de $x_kx_ix_j$ dans 
$$
\sum_{m=1}^nJ^3\left(\{x_m,x_k\}\frac{\partial \{x_i,x_j\}}{\partial x_m}\right)
$$
est 
$$
g_{i,j,E_i+E_j}\left(g_{j,k,E_k+E_j}+g_{i,k,E_k+E_i}\right).
$$
Celui de 
$$
\sum_{l=1}^pJ^3\left(X_{l,k}\frac{\partial \{x_i,x_j\}}{\partial x_{n+l}}\right)
$$
est 
$$
\sum_{l=1}^p\lambda_{l,k}\frac{\partial g_{i,j,E_i+E_j}}{\partial x_{l+n}}=\Lambda_k(g_{i,j,E_i+E_j}).
$$

On trouve alors que l'on a 
\begin{eqnarray}
\Lambda_j(g_{k,i,E_k+E_i})-\Lambda_i(g_{k,j,E_k+E_j})+\Lambda_k(g_{j,i,E_j+E_i})& = & g_{k,i,E_k+E_i}(g_{i,j,E_j+E_i}+g_{k,j,E_k+E_j})\nonumber\\
& & - g_{k,j,E_k+E_j}(g_{j,i,E_j+E_i}+g_{k,i,E_k+E_i})\nonumber\\
& & + g_{i,j,E_j+E_i}(g_{i,k,E_k+E_i}+g_{j,k,E_k+E_j})\nonumber\\
& = & 0.\label{cocycle}
\end{eqnarray}
\end{proof}
\begin{lemm}\label{resolution}
Il existe des fonctions $\gamma_1,\ldots, \gamma_n$ holomorphes au voisinage
de $K$ dans $\Bbb C^p$ et nulles \`a l'origine telles que si l'on pose 
$x_i=\exp \gamma_i(y'')y_i$, $1\leq i\leq n$, alors le coefficient de $y_iy_j$
dans $\{y_i,y_j\}'$ est constant; il est m\^eme nul si $i$ ou $j$ est inf\'erieur ou \'egal \`a $p$.
\end{lemm}
\begin{proof}
Rappelons que, pour $1\leq i\leq n$, 
$$
\Lambda_i = \sum_{k=1}^p \lambda_{k,i}\frac{\partial}{\partial x_{n+k}}.
$$
Par hypoth\`ese, on peut supposer que $\Lambda_1,\ldots,\Lambda_p$ sont
lin\'eairement ind\'ependants sur $\Bbb C$.
Quitte \`a faire un changement lin\'eaire de coordonn\'ees $y''=Bx''$, on peut supposer
que, pour $1\leq k\leq p$, 
$$
\Lambda_k = \frac{\partial}{\partial x_{n+k}}.
$$
On notera $K':=BK$ le bon compact de $0$.
Dans ces nouvelles coordonn\'ees, on a 
$$
P=\sum_{r=1}^pY_r(x',y'')\wedge \frac{\partial }{\partial y_{n+r}}+\sum_{1\leq
  i<j\leq n}\left(\sum_{\substack{Q\in \Bbb N^n_2,\\ \forall r\;
      (Q,\lambda^r)=\lambda_{r,i}+\lambda_{r,j}}}\tilde g_{i,j,Q}(y'')(x')^{Q}\right)\frac{\partial }{\partial x_{i}}\wedge \frac{\partial }{\partial x_{j}}
$$
o\`u $\tilde g_{i,j,Q}(y''):=g_{i,j,Q}(B^{-1}y'')\in {\cal O}_p(V')$, $V'=BV$,
$V$ voisinage $K$ et
$$
\begin{pmatrix}Y_1(x',y'')\\\vdots \\ Y_p(x',y'')\end{pmatrix}={} ^tB^{-1}\begin{pmatrix}X_1(x',B^{-1}y'')\\\vdots \\ Y_p(x',B^{-1}y'')\end{pmatrix}.
$$
Posons 
\begin{eqnarray*}
\gamma_{1,i}(y'') & = & \int_0^{y_{n+1}}\tilde g_{1,i,
  E_1+E_i}(t,y_{n+2},\ldots,y_N)dt \quad 1\leq i\leq n\\
x_i & = & \exp \left(\gamma_{1,i}(y'')\right) y_i \quad 1\leq i\leq n.
\end{eqnarray*}
On a $\gamma_{1,i}(0)=0$ et $\gamma_{1,1}\equiv 0$. Les $\gamma_{1,i}$ sont
holomorphes sur $(\Bbb C^n,0)\times V'$.
Si $i$ et $j$ sont distincts de $1$ et inf\'erieurs \`a $p$ , on a 
\begin{eqnarray*}
\Lambda_j(\gamma_{1,i})-\Lambda_i(\gamma_{1,j}) & = &
\int_0^{y_{n+1}}\left(\frac{\partial \tilde g_{1,i,
      E_1+E_i}(t,y_{n+2},\ldots,y_N)}{\partial x_{n+j}}-\frac{\partial \tilde g_{1,j,
      E_1+E_j}(t,y_{n+2},\ldots,y_N)}{\partial x_{n+i}}\right)dt\\
& = & \int_0^{y_{n+1}}\frac{\partial \tilde g_{i,j,
    E_i+E_j}(t,y_{n+2},\ldots,y_N)}{\partial t}dt\quad\quad(\text{gr\^ace \`a
  }(\ref{cocycle}))\\
&= & \tilde g_{i,j, E_i+E_j}(y_{n+1},y_{n+2},\ldots,y_N)-\tilde g_{i,j,
  E_i+E_j}(0,y_{n+2},\ldots,y_N).
\end{eqnarray*}
D'apr\`es $(\ref{cobord})$, le coefficient de $y_iy_j$ dans $\{y_i,y_j\}'$ est 
$$
\tilde g_{i,j, E_i+E_j}^{(2)}(y_{n+2},\ldots,y_N):= \tilde g_{i,j,
  E_i+E_j}(0,y_{n+2},\ldots,y_N)\quad\quad 2\leq i,j\leq p.
$$
De plus, si $1\leq j\leq n$, on a 
$$
\tilde g_{1,j, E_1+E_j}(y_{n+1},y_{n+2},\ldots,y_N)-\Lambda_1(\gamma_{1,j})=0.
$$
Donc, d'apr\`es $(\ref{cobord})$, le coefficient $\tilde g_{1,j,
  E_1+E_j}^{(2)}(y_{n+2},\ldots,y_N)$ de $y_1y_j$ dans $\{y_1,y_j\}'$ est nul
pour $1\leq j\leq n$. Puisque l'on transforme un crochet de Poisson en un
autre, alors les nouveaux coefficients v\'erifient encore l'\'equation de
cocycle $(\ref{cocycle})$ (les $\Lambda_i$ sont invariants par le changement
de coordonn\'ees). En y posant $k=1$, il vient 
$$
\Lambda_1\left(\tilde g_{i,j, E_i+E_j}^{(2)}\right)=0\quad\quad 1\leq i,j\leq n.
$$

On montre par r\'ecurrence sur l'entier $1\leq d\leq p$ qu'il existe des
fonctions $\gamma_{d,i}$, $1\leq i\leq n$, holomorphes sur $(\Bbb C^n,0)\times
V'$ et
nulles en $0$ telles qu'en posant
\begin{eqnarray*}
x_i & = & \exp (\gamma_{d,i}(y'')) y_i\quad\quad 1\leq i\leq n,
\end{eqnarray*}
les coefficients $g^{(d+1)}_{i,j,E_i+E_j}$ de $y_iy_j$ dans $\{y_i,y_j\}'$ v\'erifient
\begin{enumerate}
\item $g^{(d+1)}_{i,j,E_i+E_j}\equiv 0$ pour $1\leq i\leq d$, $1\leq j\leq n$,
\item les $g^{(d+1)}_{i,j,E_i+E_j}$, $1\leq i,j\leq p$, ne d\'ependent que de $x_{n+d+1},\ldots, x_{n+p}$,
\item $\Lambda_k\left(g^{(d+1)}_{i,j,E_i+E_j}\right)\equiv 0$ pour $1\leq
  k\leq d$ et $1\leq i,j\leq n$.
\end{enumerate}

D'apr\`es ce qui pr\'ec\`ede, cela est vrai pour $d=1$. Supposons que cela soit
vrai pour tout $d\leq k<p$. Posons alors 
\begin{eqnarray*}
\gamma_{k,i}(y'') & = & \int_0^{y_{n+k}}\tilde g_{k,i, E_k+E_i}^{(k)}(t,y_{n+k+1},\ldots,y_N)dt \quad 1\leq i\leq n,\\
x_i & = & \exp \left(\gamma_{k,i}(y'')\right) y_i \quad 1\leq i\leq n.
\end{eqnarray*}
Comme plus haut, on a, pour $1\leq i,j\leq p$ distincts de $k$, 
$$
\Lambda_j(\gamma_{k,i})-\Lambda_i(\gamma_{k,j})=g^{(k)}_{i,j,E_i+E_j}(y_{n+k},\ldots,y_N)-g^{(k)}_{i,j,E_i+E_j}(0,y_{n+k+1},\ldots,y_N).
$$
Le coefficient $g^{(k+1)}_{i,j,E_i+E_j}$ de $y_iy_j$ dans $\{y_i,y_j\}'$, pour
$1\leq i,j\leq p$ distincts de $k$, est donc 
$$
g^{(k+1)}_{i,j,E_i+E_j}:=g^{(k)}_{i,j,E_i+E_j}(0,y_{n+k+1},\ldots,y_N).
$$
En particulier, pour $1\leq i< k$ et $1\leq j\leq p$,
$$
g^{(k+1)}_{i,j,E_i+E_j}\equiv 0
$$
et 
$$
g^{(k+1)}_{i,j,E_i+E_j} =
g^{(k)}_{i,j,E_i+E_j}+\Lambda_i(\gamma_{k,j})-\Lambda_j(\gamma_{k,i});
$$
pour $1\leq i\leq k$ et $1\leq j\leq n$. Pour $i=k$, on obtient
$$
g^{(k+1)}_{k,j,E_k+E_j} = g^{(k)}_{k,j,E_k+E_j}-\Lambda_k(\gamma_{k,j})\equiv
0.
$$
Or, les $\gamma_{k,j}$ ne d\'ependent que de $x_{n+k},\ldots, x_{n+p}$. Donc, $\Lambda_i(\gamma_{k,j})=\partial \gamma_{k,j}/\partial x_{n+i}=0$ pour $1\leq i\leq k-1$. De plus, d'apr\`es les propri\'et\'es du crochet de Poisson, on a $\tilde g_{k,i, E_k+E_i}^{(k)}=-\tilde g_{i,k, E_k+E_i}^{(k)}$. Par hypoth\`ese de r\'ecurrence, $\tilde g_{i,k, E_k+E_i}^{(k)}$ est nul pour $i<k$; il en est donc de m\^eme pour $\gamma_{k,i}$. Par cons\'equent, pour $1\leq i \leq k$ et $1\leq j\leq n$, on a 
$$
g^{(k+1)}_{i,j,E_i+E_j} = g^{(k)}_{i,j,E_i+E_j}\equiv 0.
$$

En utilisant la relation de cocycle $(\ref{cocycle})$ dans les nouvelles
coordonn\'ees, on obtient alors 
$$
\Lambda_d\left(\tilde g_{i,j, E_i+E_j}^{(k+1)}\right)=0\quad\quad 1\leq i,j\leq n,
$$
pour $1\leq d\leq k$. Ce qui d\'emontre la r\'ecurrence.

Ainsi, pour $1\leq k\leq p$ et $1\leq i,j\leq n$, on a 
$$
\tilde g_{k,j,E_k+E_j}^{(p+1)}\equiv 0.
$$
et
$$
\Lambda_k\left(\tilde g_{i,j,E_i+E_j}^{(p+1)}\right)=0.
$$
Puisque les fonctions $\tilde g_{i,j,E_i+E_j}^{(p+1)}$ ne d\'ependent que des variables
$x_{n+1},\ldots, x_{n+p}$, alors, pour $1\leq i,j\leq n$,
$\tilde g_{i,j,E_i+E_j}^{(p+1)}$ est constant. Dans ces nouvelles
coordonn\'ees, le champ $Y_k$ est transform\'e en 
\begin{eqnarray*}
\tilde Y_k(y',y'') & = & \sum_{i=1}^n \tilde
Y_{k,i}(y',y'')\frac{\partial}{\partial
  x_i}\\
& = & \sum_{i=1}^n\exp(\gamma_i(y''))Y_{k,i}\left(\exp(-\gamma_1(y''))y_1,\ldots,\exp(-\gamma_n(y''))y_n ,y''\right)\frac{\partial}{\partial x_i}.
\end{eqnarray*}
La partie lin\'eaire reste donc inchang\'ee. Il ne reste plus qu'\`a faire
le changement de variables $x''=B^{-1}y''$.
\end{proof}

Finalement, nous avons montr\'e que, dans un bon syst\`eme de coordonn\'ees holomorphes, on peut \'ecrire

\begin{eqnarray}
P & = & \sum_{r=1}^p\left(\sum_{s=1}^p\tilde a_{r,s}(x)S_s\right)\wedge \frac{\partial }{\partial x_{n+r}}+\sum_{p+1\leq i<j\leq n}c_{i,j}x_ix_j\frac{\partial }{\partial x_{i}}\wedge \frac{\partial }{\partial x_{j}}\nonumber\\
& & +\sum_{1\leq i<j\leq n}\left(\sum_{\substack{Q\in \Bbb N^n_2,\\
Q\neq E_i+E_j\\ \forall r\;
(Q,\lambda^r)=\lambda_{r,i}+\lambda_{r,j}}}g_{i,j,Q}(x'')(x')^{Q}\right)\frac{\partial
}{\partial x_{i}}\wedge \frac{\partial }{\partial x_{j}}\label{1nf}
\end{eqnarray}
o\`u $\tilde a_{r,s}\in {\cal O}_p(K)\otimes {\cal O}_n^S$, $c_{i,j}\in \Bbb
C$ et $g_{i,j,Q}\in {\cal O}_p(V)$, $V$ \'etant un voisinage de $K$ dans $\Bbb
C^p$.

\section{Preuve du th\'eor\`eme \ref{theo2}}

Appliquons tout d'abord le th\'eor\`eme pr\'ec\'edent et posons
$$
X_i(x) := \sum_{j=1}^p \tilde a_{i,j}(x)S_j,
$$
avec $\tilde a_{i,j}\in {\cal O}_p(K)\otimes {\cal O}_n^S$. Nous d\'esignerons
par $A$ la matrice dont les \'el\'ements sont les $\tilde a_{i,j}$. Commen\c{c}ons par d\'emontrer le lemme suivant
\begin{lemm}
Il existe $p$ familles champs de vecteurs $A_i$ appartenant \`a ${\cal O}_p(K)\otimes \left(\vfo n 1\right)^S$ tel que 
$$
P = \sum_{i=1}^p X_i\wedge \left( \frac{\partial }{\partial x_{n+i}} + A_i\right)
$$
\end{lemm}
Nous poserons, dans la suite,
$$
\tilde A_i := \frac{\partial }{\partial x_{n+i}} + A_i;\quad \tilde P:= \sum_{i=1}^p X_i\wedge A_i.
$$
\begin{proof}

Nous reprenons les arguments de J.-P. Dufour \cite{dufour2}.
Par hypoth\`ese, le rang de $P$ est $2p$; donc 
$$
\underbrace{P\wedge P\wedge \ldots\wedge P}_{(p+1)\text{fois}}=0.
$$
Le terme de cette \'equation qui contient ${\partial }/{\partial x_{n+1}}\wedge\ldots\wedge{\partial }/{\partial x_{n+p}}$ donne l'\'equation
$$
\tilde P\wedge X_1\wedge\ldots\wedge X_p = 0.
$$
De plus, $X_1\wedge\ldots\wedge X_p = (\det A)\, S_1\wedge\ldots\wedge S_p$ et $\det A$ est une unit\'e dans ${\cal O}_p(K)\otimes {\cal O}_n$.
Soit $I$ l'id\'eal engendr\'e par les coefficients de $X_1\wedge\ldots\wedge X_p$ dans ${\cal O}_p(K)\otimes {\cal O}_n$. C'est l'id\'eal engendr\'e par les mon\^omes $x_{i_1}\cdots x_{i_p}$, $1\leq i_1<\cdots<i_p\leq n$.
Soient 
\begin{itemize}
\item $a:=x_1\cdots x_p$,
\item $b:=\sum_{i=1}^px_1\cdots x_{i-1}\hat x_i x_{i+1}\cdots x_{p+1}$ et
\item $c:=\sum_{1\leq i<j \leq p}^px_1\cdots x_{i-1}\hat x_i x_{i+1}\cdots x_{j-1}\hat x_j x_{j+1}\cdots x_{p+2}$.
\end{itemize}
Comme d'habitude, $\hat x_i$ signifie que la variable $x_i$ est absente; on rappel que l'on a suppos\'e $n>p+1$. Alors, $a\in I$ n'est pas un diviseur de z\'ero dans ${\cal O}_p(K)\otimes {\cal O}_n$, $b\in I$ n'est pas un diviseur de z\'ero dans ${\cal O}_p(K)\otimes {\cal O}_n/(a)$ et 
$c\in I$ n'est pas un diviseur de z\'ero dans ${\cal O}_p(K)\otimes {\cal O}_n/(a,b)$.
Par cons\'equent, la profondeur de l'id\'eal $I$ est sup\'erieure \`a $3$. 
Or, d'apr\`es le th\'eor\`eme de Frisch \cite{frisch, siu-noether}, les
anneaux ${\cal O}_p(K)$ et ${\cal O}_p(K)\otimes {\cal O}_n$ sont noeth\'eriens.
D'apr\`es le th\'eor\`eme de K. Saito \cite{saito}, il existe alors $p$ familles champs de vecteurs $A_i$ appartenant \`a ${\cal O}_p(K)\otimes \vf n $ tel que 
$$
\tilde P = \sum_{i=1}^p X_i\wedge A_i.
$$
Puisque $\tilde P$ est d'ordre sup\'erieur \`a $2$ en $x'$, on peut choisir
$A_i$ nul en $(0,x'')$. Montrons alors que l'on peut choisir des $A_i$ qui
commutent avec les $S_i$. En effet, on a
\begin{eqnarray*}
\tilde P  & = &\sum_{p+1\leq i<j\leq n}c_{i,j}x_ix_j\frac{\partial }{\partial x_{i}}\wedge \frac{\partial }{\partial x_{j}} \\
& & +\sum_{1\leq i<j\leq n}\left(\sum_{\substack{Q\in \Bbb N^n_2,\\
Q\neq E_i+E_j\\ \forall r\;
(Q,\lambda^r)=\lambda_{r,i}+\lambda_{r,j}}}\tilde g_{i,j,Q}(X'')(X')^{Q}\right)\frac{\partial
}{\partial x_{i}}\wedge \frac{\partial }{\partial x_{j}}.
\end{eqnarray*}
Par cons\'equent, pour tout $g\in {\lie g}$, $[S(g),\tilde P]=0$ (d'apr\`es  le lemme \ref{poids-poisson}). On a aussi
$$
[S(g),\tilde P]= \sum_{i=1}^p X_i\wedge [S(g),A_i],
$$
car les $X_i$ commutent aux $S_j$. D\'ecomposons ensuite chaque $A_i$ selon les espaces de poids de $S$ dans l'espaces des champs de vecteurs (en fait, il faudrait faire cette d\'ecomposition dans chaque espace de champs de bi-vecteurs de $\Bbb C^n$ homog\`enes \`a coefficients dans ${\cal O}_p(K)$). On a alors 
\begin{eqnarray*}
0=[S(g),\tilde P] & = &\sum_{i=1}^p X_i\wedge \sum_{\alpha}\alpha(g)A_{i,\alpha}\\
& = & \sum_{\alpha}\alpha(g)\left(\sum_{i=1}^p X_i\wedge A_{i,\alpha}\right).
\end{eqnarray*}
Puisque les $X_i$ sont de poids nuls, l'\'egalit\'e pr\'ec\'edente d\'efinit la d\'ecomposition de $[S(g),\tilde P]$ selon les espaces de poids de $S$ dans l'espaces des germes de champs de bi-vecteurs nul \`a l'origine. Les espaces de cette d\'ecomposition \'etant en somme directe, on a, pour tout $\alpha\neq 0$, 
$$
\sum_{i=1}^p X_i\wedge A_{i,\alpha}=0,
$$
et on a 
$$
\tilde P = \sum_{i=1}^p X_i\wedge A_{i}=\sum_{i=1}^p X_i\wedge A_{i,0}.
$$
Par d\'efinition, pour tout $g\in {\lie g}$, $[S(g),A_{i,0}]=0$.
\end{proof}

Puisque $P$ est une structure de Poisson, on a $[P,P]=0$. En d\'eveloppant, on obtient
\begin{eqnarray}
0=\sum_{i,j=1}^p [X_i\wedge\tilde A_i, X_j\wedge\tilde A_j] & = & \sum_{i,j=1}^p [X_i\wedge\tilde A_i, X_j]\wedge \tilde A_j-X_j\wedge[X_i\wedge\tilde A_i, \tilde A_j]\nonumber\\
& = & \sum_{i,j=1}^p [X_j,X_i]\wedge \tilde A_i\wedge \tilde A_j + X_i\wedge [X_j, \tilde A_i]\wedge \tilde A_j\nonumber\\
& & + \sum_{i,j=1}^p -X_j\wedge[X_i,\tilde A_j]\wedge \tilde A_i - X_j\wedge X_i\wedge[\tilde A_j, \tilde A_i]\nonumber\\
& = & 2\left(\sum_{i,j=1}^p X_i\wedge [X_j, \tilde A_i]\wedge \tilde
  A_j\right.\nonumber\\
& & -\left.\sum_{1\leq i<j\leq p} X_j\wedge X_i\wedge[\tilde A_j, \tilde
  A_i]\right).\label{poisson2}
\end{eqnarray}
Prenons le terme qui contient $\partial/\partial x_{n+j}$, il verifie alors
\begin{eqnarray}\label{poisson1}
\sum_{i=1}^p X_i\wedge [X_j, \tilde A_i]&= &0.
\end{eqnarray}
On en d\'eduit
\begin{equation}\label{poisson3}
\sum_{1\leq i<j\leq p} X_j\wedge X_i\wedge[\tilde A_j, \tilde
  A_i]= 0.
\end{equation}
En faisant le produit de $(\ref{poisson1})$ par $X_1\wedge \ldots \wedge X_{j-1}\wedge X_{j+1}\wedge\ldots\wedge  X_p$, on obtient 
$$
[X_j, \tilde A_i]\wedge X_1\wedge \ldots \wedge X_p=0.
$$
En appliquant le th\'eor\`eme de Saito, il existe alors des famille de fonctions $\theta_{i,j}^l\in {\cal O}_p(K)\otimes {\cal O}_n$ telles que
\begin{eqnarray}\label{theta}
[\tilde A_i,X_j]=\sum_{l=1}^p \theta_{i,j}^l X_l.
\end{eqnarray}
De surcro\^{\i}t, on a 
$$
[S(g),[\tilde A_i,X_j]]  = -[\tilde A_i, [X_j, S(g)]] - [X_j,[S(g),\tilde A_i]].
$$
Or, par d\'efinition, $[X_j,S(g)]=0$ et par le lemme pr\'ec\'edent, $[S(g),\tilde A_i]=0$ ($S(g)$ ne d\'epend pas de $x''$). On en d\'eduit que $S(g)(\theta_{i,j}^l)=0$; c'est-\`a-dire $\theta_{i,j}^l\in {\cal O}_p(K)\otimes {\cal O}_n^S$. En reinjectant l'expression $(\ref{theta})$ dans l'\'equation $(\ref{poisson1})$, on obtient que $\theta_{i,j}^l=\theta_{l,j}^i$. Prenons le $1$-jet de l'\'equation $(\ref{theta})$ et posons $J^1(A_i)=\sum_{k,l=1}^n A_{i,k,l}^1(x'')x_k\partial/\partial x_l$. On a alors
$$
\sum_{k,l=1}^n(\lambda_{j,k}-\lambda_{j,l})A_{i,k,l}^1(x'')x_k\frac{\partial}{\partial x_l}=\sum_{k=1}^n\left(\sum_{l=1}^p \theta_{i,j}^l(0,x'')\lambda_{l,k}\right)x_k \frac{\partial}{\partial x_k}.
$$
Donc, pour $1\leq k\leq n$, on a
$$
\sum_{l=1}^p \theta_{i,j}^l(0,x'')\lambda_{l,k}=0.
$$
Or, par hypoth\`ese, la famille $\{S_i\}_{i=1}^p$ est lin\'eairement ind\'ependante sur $\Bbb C$. Donc, pour tout $i,j,l$,  $\theta_{i,j}^l(0,x'')\equiv 0$.

En faisant le produit de $(\ref{poisson3})$ par $X_1\wedge\ldots \wedge X_{i-1}\wedge X_{i+1}\wedge
\ldots \wedge X_{j-1}\wedge X_{j+1}\wedge\ldots\wedge X_p$, on obtient
$$
[\tilde A_j, \tilde A_i]\wedge X_1\wedge \ldots \wedge X_p=0.
$$
En appliquant le th\'eor\`eme de Saito, il existe alors des familles de fonctions $\gamma_{i,j}^l\in {\cal O}_p(K)\otimes {\cal O}_n$ telles que
\begin{eqnarray}\label{gamma}
[\tilde A_i,\tilde A_j]=\sum_{l=1}^p \gamma_{i,j}^l X_l.
\end{eqnarray}
De surcro\^{\i}t, on a 
$$
[S_k,[\tilde A_i,\tilde A_j]]  = -[\tilde A_i, [\tilde A_j, S_k]] - [\tilde A_j,[S_k,\tilde A_i]]=0.
$$
On en d\'eduit que, pour $1\leq k\leq p$, $S_k(\gamma_{i,j}^l)=0$; donc
$\gamma_{i,j}^l\in {\cal O}_p(K)\otimes {\cal O}_n^S$. 

Reinjectons l'\'egalit\'e $(\ref{gamma})$ dans l'\'equation
$(\ref{poisson2})$. On obtient 
$$
\sum_{1\leq i<j\leq p}\sum_{l=1}^p \gamma_{i,j}^l X_l\wedge X_i\wedge X_j = 0.
$$
En prenant les coefficients de $ X_l\wedge X_i\wedge X_j $ avec $l<i<j$, on en d\'eduit que l'on a
\begin{equation}
\gamma_{i,j}^l=\gamma_{l,j}^i.\label{gamma-sym}
\end{equation}
%

Montrons par r\'ecurrence sur l'entier $1\leq q\leq p$  qu'il existe un
changement de coordonn\'ees holomorphes dans lequel $\tilde A_i=\partial/
\partial x_{n+i}$, pour $1\leq i\leq q$. 

Pour $q=1$, consid\'erons le flot $\phi_{1,t}(x',x'')$ de $\tilde A_1$ et posons 
$$
\psi_1(x',x''):=\phi_{1,x_{n+1}}(x',0,x_{n+2},\ldots, x_N). 
$$
C'est un
diff\'eomorphisme local au voisinage de l'origine qui est tangent \`a
l'identit\'e en ce point et qui v\'erifie $(\psi_1^{-1})_*\tilde A_1= \partial/\partial x_{n+1}$. Ce flot v\'erifie le syst\`eme d'\'equations diff\'erentielles suivantes
\begin{eqnarray*}
\frac{d \psi_{1,t,i}}{dt}(x',x'') & = & A_i(\psi_{1,t}(x))\quad i=1,\ldots, n\\
\frac{d \psi_{1,t,n+1}}{dt}(x',x'') & = & 1\\
\frac{d \psi_{1,t,n+j}}{dt}(x',x'') & = & 0\quad j=2,\ldots, p,
\end{eqnarray*}
o\`u $\psi_{1,t,i}$ d\'esigne la $i$\`eme coordonn\'ees de $\psi_{1,t}$. En
particulier, $\psi_1$ laisse invariante chacune des coordonn\'ees $x_{n+j}$,
$j=1,\ldots, p$. Appliquons la m\'ethode du chemin
\cite{Chaperon-ast}[apendice 0]; on obtient, pour $1\leq k\leq p$,
$$
\frac{d ((\psi^{-1}_1)_*S_k)}{d x_{n+1}}= (\psi^{-1}_1)_*\left([\tilde A_1, S_k]\right)= 0.
$$
On en d\'eduit que
$$
(\psi^{-1}_1)_*S_k = (S_k)_{|x_{n+1}=0}= S_k;
$$
c'est-\`a-dire, $\psi^{-1}_1$ pr\'eserve chacun des champs $S_k$. Soit alors $f\in {\cal O}_p(K)\otimes {\cal O}_n^S$ alors $(\psi_1^{-1})_*f=f\circ \psi^{-1}_1\in{\cal O}_p(K)\otimes {\cal O}_n^S$. En effet, 
$$
0=(\psi^{-1}_1)_*(S(g)(f))=(\psi^{-1}_1)_*S(g)(f\circ \psi^{-1}_1)=S(g)(f\circ\psi^{-1}_1).
$$
Il en est de m\^eme pour les \'el\'ements de ${\cal O}(K)\otimes \left( \vfo
n 1\right)^S$. On a alors 
$$
X'_k(x):=(\psi^{-1}_1)_*X_k = \sum_{l=1}^p (\tilde a_{k,l}\circ \psi^{-1}_1) S_l.
$$
On posera $\tilde a_{k,l}'=\tilde a_{k,l}\circ \psi^{-1}_1$.
Par cons\'equent, dans ces nouvelles coordonn\'ees, on peut \'ecrire
$$
P= \sum_{i=1}^p X_i'(x)\wedge \tilde A_i'
$$
avec $\tilde A_1'=\partial/\partial x_{n+1}$ et $[\tilde A_1', S_k']=0$. Afin d'all\'eger le texte, nous continuerons \`a \'ecrire $X_i$ \`a la place de $X_i'$ (etc...) lorsque cela ne pr\^etera pas \`a confusion.

Supposons que l'assertion soit vraie pour $1\leq i\leq q-1$, c'est-\`a-dire
 $A_1=\cdots=A_{q-1}=0$.

Commen\c{c}ons par faire un changement de coordonn\'ees qui transforme $\tilde
A_q$ en $\partial/\partial x_{n+q}+ \sum_{j=1}^p c_{q,j}X_j$ o\`u $c_{i,j}\in {\cal O}_p(K')\otimes {\cal O}_n^S$, $K'$ \'etant un
bon compact de l'origine de $\Bbb C^p$. On demande aussi que ce changement
n'affecte pas ni les $\tilde A_i=\partial/\partial x_{n+i}$, $1\leq i\leq q-1$, ni les
$S_j$, $1\leq j\leq p$. Pour ce faire, cherchons des fonctions $\beta_{q,l}$ telles
que, pour $1\leq i\leq q-1$, 
$$
0=\left [\tilde A_i, \tilde A_q -\sum_{j=1}^p \beta_{q,j} X_j\right ]= \sum_{j=1}^p
\left (\gamma_{i,q}^j-\frac{\partial \beta_{q,j}}{\partial x_{n+i}}\right )X_j-\beta_{q,j}[\tilde A_i, X_l].
$$
Nous avons utilis\'e les \'equations $(\ref{gamma})$ dans la derni\`ere \'egalit\'e. En utilisant les \'egalit\'es $(\ref{theta})$; 
on obtient, pour $1\leq k\leq p$ et $1\leq i\leq q-1$,
\begin{equation}\label{betai}
\frac{\partial \beta_{q,k}}{\partial x_{n+i}}  =  - \left(\sum_{l=1}^p\theta_{i,l}^k\beta_{q,l}\right)+\gamma_{i,q}^k.
\end{equation}
On \'ecrira ces \'equations sous forme matricielle de la mani\`ere suivante
\begin{equation}\label{sol-beta}
\frac{\partial \beta_{q}}{\partial x_{n+i}}  =  - \Theta_i\beta_{q}+\gamma_{i,q}\quad\quad i=1,\ldots, q-1
\end{equation}
o\`u $\beta_q$ d\'esigne le vecteur de coordonn\'ees $\beta_{q,1},\ldots , \beta_{q, p}$, $\Theta_i=\left(\theta_{i,l}^k\right)_{1\leq k,l\leq p}$ et 
$\gamma_{i,q}$ d\'esigne le vecteur de coordonn\'ees $\gamma_{i,q}^1,\ldots,\gamma_{i,q}^p$.
Pour $1\leq i,j <q$, $\tilde A_i$ et $\tilde A_j$ commutent entre eux. Par
l'\'egalit\'e de Jacobi, on a alors, pour $1\leq i,j <q$,
\begin{eqnarray*}
[\tilde A_i,[\tilde A_j,\tilde A_q]]&=& [\tilde A_j,[\tilde A_i,\tilde A_q]],
\end{eqnarray*}
soit,
\begin{eqnarray*}
\left[\tilde A_i, \sum_{k=1}^p \gamma_{j,q}^k X_k\right]& = & \left[\tilde A_j, \sum_{k=1}^p \gamma_{i,q}^k X_k\right]\quad\text{(par $(\ref{gamma})$)},
\end{eqnarray*}
c'est-\`a-dire, 
\begin{eqnarray*}
\sum_{k=1}^p \frac{\partial \gamma_{j,q}^k}{\partial x_{n+i}}X_k+\gamma_{j,q}^k[\tilde A_i, X_k] &= & \sum_{k=1}^p \frac{\partial \gamma_{i,q}^k}{\partial x_{n+j}}X_k+\gamma_{i,q}^k[\tilde A_j, X_k],
\end{eqnarray*}
soit, finalement, 
\begin{eqnarray*}
\sum_{k=1}^p \left(\frac{\partial \gamma_{j,q}^k}{\partial x_{n+i}}+\left(\sum_{l=1}^p\gamma_{j,q}^l\theta_{i,l}^k\right)\right)X_k &= & \sum_{k=1}^p \left(\frac{\partial \gamma_{i,q}^k}{\partial x_{n+j}}+\left(\sum_{l=1}^p\gamma_{i,q}^l\theta_{j,l}^k\right)\right)X_k \quad\text{(par $(\ref{theta})$)}.
\end{eqnarray*}
En d'autres termes, les \'equations suivantes sont satisfaites
\begin{equation}\label{compat1}
\frac{\partial \gamma_{i,q}}{x_{n+j}}+\Theta_j\gamma_{i,q}=\frac{\partial
  \gamma_{j,q}}{x_{n+i}}+\Theta_i\gamma_{j,q}, \quad 1\leq i,j <q.
\end{equation}

D'autre part, par l'identit\'e de Jacobi, on a aussi, pour tous $1\leq i,j\leq q-1$ et $1\leq l\leq p$,
\begin{eqnarray*}
[\tilde A_i,[\tilde A_j,X_l]]&=& [\tilde A_j,[\tilde A_i,X_l]],
\end{eqnarray*}
soit,
\begin{eqnarray*}
\left[\tilde A_i, \sum_{k=1}^p \theta_{j,l}^k X_k\right]& = & \left[\tilde A_j, \sum_{k=1}^p \gamma_{i,l}^k X_k\right]\quad\text{(par $(\ref{theta})$)},
\end{eqnarray*}
c'est-\`a-dire, 
\begin{eqnarray*}
\sum_{k=1}^p \left(\frac{\partial \theta_{j,l}^k}{\partial x_{n+i}}+\left(\sum_{m=1}^p\theta_{j,l}^m\theta_{i,m}^k\right)\right)X_k &= & \sum_{k=1}^p \left(\frac{\partial \theta_{i,l}^k}{\partial x_{n+j}}+\left(\sum_{m=1}^p\theta_{i,l}^m\theta_{j,m}^k\right)\right)X_k \quad\text{(par $(\ref{theta})$)}.
\end{eqnarray*}
En d'autres termes, les \'equations suivantes sont satisfaites
\begin{equation}\label{compat2}
\frac{\partial \Theta_{j}}{x_{n+i}}+\Theta_i\Theta_j=\frac{\partial
  \Theta_{i}}{x_{n+j}}+\Theta_j\Theta_i, \quad 1\leq i,j <q.
\end{equation}

En utilisant les \'equations $(\ref{compat1})$ et $(\ref{compat2})$, on obtient, pour $1\leq i,j<q$ et pour tout $v\in \Bbb C^p$,
$$
-\frac{\partial \Theta_{j}}{x_{n+i}}v+ \frac{\partial \gamma_{j,q}}{\partial x_{n+i}}+\frac{\partial \Theta_{i}}{x_{n+j}}v -\frac{\partial \gamma_{i,q}}{\partial x_{n+j}}+\Theta_i(-\Theta_jv+ \gamma_{j,q})-\Theta_j(-\Theta_iv+ \gamma_{i,q})=0
$$

Ces \'equations ne sont autres que les {\bf \'equations de compatibilit\'e} du syst\`eme $(\ref{sol-beta})$ (quitte \`a rajouter des \'equations "vides" $\partial \beta_q/\partial z_r=0$, $r=1,\ldots, p-q+1$ afin de rendre "carr\'e" le syst\`eme).
Par cons\'equent, le syst\`eme d'\'equations $(\ref{sol-beta})$ admet une unique solution holomorphe sur un voisinage $U'\subset K$ de l'origine de $\Bbb C^p$ et nulle en ce point \cite{dieudonne1}[th\'eor\`eme de Frobenius 10.9.4]. Ce voisinage contient un bon compact $K'$ de $0$. De plus, cette solution d\'epend holomorphiquement des param\`etres $x'$ dans un certain voisinage de l'origine de $\Bbb C^n$. D'apr\`es ce qui pr\'ec\`ede, les $\gamma_{i,q}^j$ (resp. $\theta_{i,k}^l$) appartiennent \`a ${\cal O}_p(K)\otimes {\cal O}_n^S$. 
Alors les $\beta_{q,j}$ appartiennent \`a ${\cal O}_p(K)\otimes {\cal O}_n^S$. En effet, les op\'erateurs de d\'erivation $\partial/\partial x_{n+i}$ et $S_k$ commutent deux \`a deux. En prenant la d\'eriv\'ee de Lie le long de $S_k$ des \'equations $(\ref{sol-beta})$, on obtient
$$
\frac{\partial S_k(\beta_{q})}{\partial x_{n+i}}  =  - \Theta_iS_k(\beta_{q})\quad\quad i=1,\ldots, q-1,
$$
syst\`eme qui admet une unique solution nulle \`a l'origine-\`a savoir- $S_k(\beta_{q})\equiv 0$.

Ainsi, on a, pour $1\leq i\leq p$, 
$$
\left[S_i, \tilde A_q -\sum_{j=1}^p \beta_{q,j} X_j\right ]= 0.
$$
On peut donc redresser, par un diff\'eomorphisme holomorphe $G_q$, $\tilde A_q
-\sum_{j=1}^p \beta_{q,j} X_j$ en $\partial/\partial x_{n+q}$ tout en
pr\'eservant les $\partial/\partial x_{n+i}$, $1\leq i\leq q-1$, les $S_j$, $1\leq j\leq p$ ainsi que les coordonn\'ees $x_{n+i}$, $\leq i\leq p$. Dans ces nouvelles coordonn\'ees, $\tilde A_q':=(G_q)_*\tilde A_q$ s'\'ecrit 
\begin{equation}\label{aq}
\frac{\partial}{\partial x_{n+q}}+ \sum_{j=1}^p c_{q,j}X_j',
\end{equation}
avec $c_{q,j}=\beta_{q,j}\circ G_q \in {\cal O}_p(K')\otimes {\cal O}_n^S$,
$X_i':=(G_q)_*X_i=\sum_{l=1}^p\tilde a_{i,l}'S_l$ et $\tilde a_{i,l}'=\tilde a_{i,l}'\circ G_q \in {\cal O}_p(K')\otimes {\cal O}_n^S$.

On ne modifie pas $P$ si l'on rajoute $\sum_{l=1}^pf_{i,l}X_l'$ \`a $\tilde A_i'$, $1\leq i\leq p$, lorsque $f_{i,l}=f_{l,i}$ : 
$$
\sum_{i=1}^p X_i'\wedge\left( \tilde A_i'+\sum_{l=1}^pf_{i,l}X_l'\right)= \sum_{i=1}^p X_i'\wedge\tilde A_i'.
$$
Cherchons alors de telles fonctions v\'erifiant en outre
$$
\left[\tilde A_i'+\sum_{l=1}^pf_{i,l}X_l',\tilde A_j'+\sum_{l=1}^pf_{j,l}X_l'\right]=0
$$
pour $1\leq i,j\leq q$ et $f_{i,l}=0$ si $i>q$. En utilisant
$(\ref{theta})$, on obtient pour tous $1\leq i,j<q$,
$$
\frac{\partial f_{j,k}}{\partial x_{n+i}}+ \sum_{l=1}^p f_{j,l}\theta_{i,l}^k-\frac{\partial f_{i,k}}{\partial x_{n+j}}+ \sum_{l=1}^p f_{i,l}\theta_{j,l}^k=0;
$$
c'est-\`a-dire, en posant $f_i$ le vecteur de coorodonn\'ees $f_{i,j}$,  
\begin{equation}\label{commuti}
\frac{\partial f_{j}}{\partial x_{n+i}}-\Theta_jf_i-\frac{\partial f_{i}}{\partial x_{n+j}}+\Theta_if_j=0.
\end{equation}
Lorsque $1\leq i<q$ et $j=q$, on obtient 
\begin{equation}\label{commutq}
\frac{\partial f_{q}}{\partial x_{n+i}}-\Theta_qf_i-\frac{\partial f_{i}}{\partial
  x_{n+q}}+\Theta_if_q=-\Theta_i c_q-\frac{\partial c_{q}}{\partial x_{n+i}},
\end{equation}
car, gr\^ace \`a $(\ref{aq})$, il suffit de rajouter $c_{q,i}$ \`a $f_{q,i}$
dans l'\'equation $(\ref{commuti})$ pour obtenir le r\'esultat.

Soit $g_q$ le vecteur de coordonn\'ees $g_{q,i}\in {\cal O}_p(K')\otimes
{\cal O}_n^S$ et posons
\begin{equation}\label{gq}
\left\{\begin{array}{rcl}\frac{\partial g_{q}}{\partial x_{n+i}}-\Theta_ig_q &  = & f_i\quad\quad
    1\leq i<q\\
\frac{\partial g_{q}}{\partial x_{n+q}}-\Theta_qg_q & = & f_q+c_q\end{array}\right.
\end{equation}
Alors, les $f_i$ sont solutions du syst\`eme d'\'equations
$(\ref{commuti})-(\ref{commutq})$. En effet, ces derni\`eres ne sont autres
que les \'equations de compatibilit\'es du syst\`eme $(\ref{gq})$ d'inconnue $g_q$.

Il nous reste \`a montrer qu'il existe $g_q$ tel que les $f_{i,j}$ v\'erifient la
condition de symm\'etrie $f_{i,j} =f_{j,i}$. Cela s'\'ecrit, pour $1\leq i,j<q $
$$
\frac{\partial g_{q,j}}{\partial x_{n+i}}-\sum_{k=1}^p\theta_{i,k}^jg_{q,k} = \frac{\partial g_{q,i}}{\partial x_{n+j}}-\sum_{k=1}^p\theta_{j,k}^ig_{q,k}.
$$
Pour $1\leq i<q$, on obtient
$$
\frac{\partial g_{q,q}}{\partial x_{n+i}}-\sum_{k=1}^p\theta_{i,k}^qg_{q,k} = \frac{\partial g_{q,i}}{\partial x_{n+q}}-\sum_{k=1}^p\theta_{q,k}^ig_{q,k}-c_{q,i}.
$$
Or, comme on l'a vu plus haut, $\theta_{i,k}^q=\theta_{q,k}^i$. Les
\'equations pr\'ec\'edentes sont donc \'equivalentes au syst\`eme 
\begin{equation}\label{gq2}
\left\{\begin{array}{rcl}\frac{\partial g_{q,j}}{\partial x_{n+i}} &  = & \frac{\partial g_{q,i}}{\partial x_{n+j}}\quad\quad
    1\leq i,j<q\\
\frac{\partial g_{q,q}}{\partial x_{n+i}} & = & \frac{\partial g_{q,i}}{\partial x_{n+q}}-c_{q,i}\end{array}\right.
\end{equation}

Posons
\begin{eqnarray*}
d_{i,l} & := & 0\quad\quad 1\leq i,l<q,\\
d_{q,l} & := & c_{q,l}\quad l<q.
\end{eqnarray*}

Si $i,j,k<q$, on a 
$$
\frac{\partial d_{j,q}}{\partial x_{n+i}}-\frac{\partial d_{k,i}}{\partial x_{n+j}}+\frac{\partial d_{i,j}}{\partial x_{n+k}}=0,
$$
car $d_{j,k}=d_{i,j}=d_{k,i}=0$.
Lorsque $k= q$, on obtient 
$$
\frac{\partial d_{j,q}}{\partial x_{n+i}}-\frac{\partial d_{q,i}}{\partial x_{n+j}}+\frac{\partial d_{i,j}}{\partial x_{n+q}}=\frac{\partial c_{q,i}}{\partial x_{n+j}}-\frac{\partial c_{q,j}}{\partial x_{n+i}}.
$$

Comme nous l'avons vu dans la preuve du lemme \ref{resolution}, pour
r\'esoudre le syst\`eme $(\ref{gq2})$ d'inconnues $g_{q,i}$, il suffit de
v\'erifier les conditions de compatibilit\'es 
\begin{equation}\label{compatgq2}
\frac{\partial c_{q,j}}{\partial x_{n+i}}= \frac{\partial c_{q,i}}{\partial
  x_{n+j}}\quad i,j< q.
\end{equation}
Or, lorsque $1\leq i< q$, on a, d'une part, 
$$
[\tilde A_i',\tilde A_q']= \sum_{k=1}^p\gamma_{i,q}^kX_k'\quad\text{(d'apr\`es ($\ref{gamma}$))};
$$
et d'autre part 
\begin{eqnarray*}
[\tilde A_i',\tilde A_q'] & =& \left[\tilde A_i',\frac{\partial}{\partial
    x_{n+q}}+\sum_{l=1}c_{q,l}X_l'\right]\\
& = & \sum_{k=1}^p \left(\frac{\partial c_{q,k}}{\partial
    x_{n+i}}+\sum_{l=1}^p \theta_{i,l}^kc_{q,l}\right)X_k'.
\end{eqnarray*}
Or, d'apr\`es $(\ref{gamma-sym})$, $\gamma_{i,q}^k=\gamma_{k,q}^i$ lorsque
$i,k<q$. Puisque $\theta_{i,k}^q=\theta_{q,k}^i$ pour tous les $i,k$ alors les
\'equations de compatibilit\'e $(\ref{compatgq2})$ sont satisfaites. Par
cons\'equent, le syst\`eme $(\ref{gq2})$ admet des solutions holomorphes
$g_{i,j}$ au voisinage $V$ de l'origine de $\Bbb C^p$. En prenant un developpement de
Taylor au voisinage de l'origne de $\Bbb C^n$, on obtient que chaque $g_{q,j}$
appartient \`a ${\cal O}_p(K')\otimes {\cal O}_n ^S$. Les fonctions $f_i$ sont alors donn\'ees par les \'equations $(\ref{gq})$ et appartiennent \`a ${\cal O}_p(K')\otimes {\cal O}_n ^S$.

Par cons\'equent, on peut redresser simultan\'ement, par un diff\'eomorphisme $\Phi_q$, les champs $\tilde
A_i'+\sum_{l=1}^pf_{i,l}X_l'$, $1\leq i\leq q$ en $\partial/\partial
x_{n+i}$. De plus, puisque chaque $\tilde
A_i'+\sum_{l=1}^pf_{i,l}X_l'$ commute aves les $S_j$ alors $\Phi_q$ pr\'eserve
les $S_j$. Enfin, il pr\'eserve chaque coordonn\'ees $x_{n+i}$, $1\leq i\leq p$. Par cons\'equent, 
$$
(\Phi_q)_*X_i' =\sum_{j=1}^p a_{i,j}''S_j\text{ avec }\;a_{i,j}''\in {\cal
  O}_p(V)\otimes {\cal O}_n ^S.
$$
Ceci d\'emontre la r\'eccurence.

On a donc montr\'e qu'il existe un bon compact $\tilde K$ de l'origine de $\Bbb C^p$ et
un diff\'eomorphisme holomorphe $\Psi$ de $(\Bbb C^N,0)$ fixant l'origine et
tangent \`a l'identit\'e en ce point tel que
$$
\Psi_*P=\sum_{i=1}^p\left(\sum_{j=1}^p a_{i,j}(x)S_j\right)\wedge \frac{\partial}{\partial
    x_{n+i}}\quad\quad a_{i,j}\in {\cal
  O}_p(\tilde K)\otimes {\cal O}_n ^S.
$$
Il reste \`a montrer qu'en faisant un changement de variables {\it ad-hoc}, on
peut supprimer la d\'ependance en $x''$ des $a_{i,j}$.

Posons $Y_i:=\sum_{j=1}^p a_{i,j}(x)S_j$. $\Psi_*P$ est une structure de
Poisson donc $[\Psi_*P,\Psi_*P]=0$. Un calcul simple montre que l'on a 
$$
[\Psi_*P,\Psi_*P]=2\sum_{i=1}^p\frac{\partial}{\partial x_{n+i}}\wedge\left(\sum_{j=1}^p\frac{\partial Y_i}{\partial x_{n+j}}\wedge Y_j\right).
$$
Par cons\'equent, $\Psi_*P$ est une structure de Poisson si et seulement si, pour $1\leq i\leq p$, 
$$
\sum_{j=1}^p\frac{\partial Y_i}{\partial x_{n+j}}\wedge Y_j= \sum_{1\leq r< s\leq p}\left(\sum_{j=1}^p\frac{\partial a_{i,r}}{\partial x_{n+j}} a_{j,s}-\frac{\partial a_{i,s}}{\partial x_{n+j}} a_{j,r}\right)S_r\wedge S_s=0.
$$
Posons alors, pour $1\leq k\leq p$,
$$
\Gamma_k=\sum_{i=1}^pa_{i,k}\frac{\partial }{\partial x_{n+i}}.
$$
Les \'equations pr\'ec\'edentes se lisent alors: $[\Gamma_r,\Gamma_s]=0$ pour $1\leq r,s\leq p$.

Montrons, dans cette situation, qu'il existe un syst\`eme de coordonn\'ees holomorphes dans lequel les $Y_i$ deviennent "ind\'ependants" de $x''$. 

En effet, la famille $\{\Gamma_k(0)\}_{1\leq k\leq l}$ est libre. Or, on a $\widehat {\cal O}_n^S=\Bbb C[[(x')^{R_1},\ldots, (x')^{R_t}]]$. Ainsi, on a $a_{i,k}=b_{i,k}((x')^{R_1},\ldots, (x')^{R_t}, x'')$ o\`u $b_{i,k}(w,x'')$ est holomorphe dans un voisinage de l'origine dans $\Bbb C^t\times \tilde K$.
Il existe un voisinage de l'origine $W$ dans $\Bbb C^t$ tel que, pour tout $w\in W$, la famille  $\{\Gamma_k(w,0)\}_{1\leq k\leq l}$ soit libre (par abus de language, nous avons \'ecrit $\Gamma_k(x)=\Gamma_k((x')^{R_1},\ldots, (x')^{R_t},x'')$).
Il existe alors un changement de variables holomorphes $X''=\psi(w,x'')$ d\'ependant holomorphiquement de $w\in W$ tel que $\psi_*\Gamma_k(w,X'')=\Gamma_k(w,0)$, pour $1\leq k\leq p$. D'autre part, les champs $\Gamma_k$ commutent aux champs $S_k$. Le diff\'eomorphisme $\psi$ laisse donc invariant chacun des $S_k$ (par la m\'ethode du chemin). Posons alors 
$(X',X'')=\Psi(x)=(x',\psi((x')^{R_1},\ldots, (x')^{R_t},x''))$. On obtient alors 
\begin{eqnarray*}
\Psi_*\left(\sum_{i=1}^p Y_i\wedge \frac{\partial}{\partial x_{n+i}}\right)(X) & = & \Psi_*\left(\sum_{i=1}^p S_i\wedge \Gamma_i\right)\\
& = & \sum_{i=1}^p S_i\wedge \Gamma_i((X')^{R_1},\ldots, (X')^{R_t},0)\\
& = & \sum_{i=1}^p \left(\sum_{j=1}^p b_{i,j}((X')^{R_1},\ldots, (X')^{R_t},0)S_j\right)\wedge \frac{\partial}{\partial X_{n+i}}.
\end{eqnarray*}
Ceci montre ainsi le r\'esultat annonc\'e.





\bibliographystyle{alpha}
\bibliography{normal,math,asympt,analyse,stolo,lie,poisson}

\end{document}